\documentclass[]{article}

\usepackage{amsfonts}
\usepackage{amsmath}
\usepackage{amssymb}
\usepackage{amsthm}
\usepackage{graphicx}
\usepackage{booktabs}
\usepackage{indentfirst}
\usepackage{placeins}
\usepackage{needspace}
\usepackage[hidelinks]{hyperref}
\usepackage{url}

\graphicspath{{figures/}}

\numberwithin{equation}{section}
\theoremstyle{plain}
\newtheorem{theorem}{Theorem}[section]
\newtheorem{proposition}[theorem]{Proposition}

\theoremstyle{remark}
\newtheorem{remark}[theorem]{Remark}

\newcommand{\R}{\mathbb{R}}
\newcommand{\one}{\mathbf 1}
\newcommand{\KL}{\operatorname{KL}}
\newcommand{\rank}{\operatorname{rank}}
\newcommand{\range}{\operatorname{range}}
\newcommand{\diag}{\operatorname{diag}}
\newcommand{\spn}{\operatorname{span}}

\newcommand{\DynamicRuleError}{4.2\times10^{-4}}

\newcommand{\MatchedDynamicSpeedMinPrecise}{1.36}
\newcommand{\MatchedDynamicSpeedMaxPrecise}{2.89}
\newcommand{\MatchedGridMax}{1024}
\newcommand{\MatchedFineGridDifference}{0.0063}
\newcommand{\MatchedNeuralFineF}{0.328}
\newcommand{\MatchedDynamicFineF}{0.240}

\newcommand{\FinalECMInvariant}{6.00\times10^{-14}}
\newcommand{\FinalECMEntropyIdentity}{3.12\times10^{-13}}
\newcommand{\FinalECMIterationMin}{2.42}
\newcommand{\FinalECMIterationMax}{2.74}

\hypersetup{
 pdftitle={Data-Adaptive Entropy-Based Moment Reduction for
 Vlasov--Poisson--Dougherty Dynamics},
 pdfauthor={Shiheng Zhang, Jiahao Zhang, Guang Lin},
 pdfkeywords={Vlasov--Poisson--Dougherty system; entropy-based moment
 method; structure-preserving discretization; reduced-order model;
 hyper-reduction; proper orthogonal decomposition}
}

\newcommand{\DomainFullLongDelta}{6.42\times10^{-12}}

\begin{document}

\title{Data-Adaptive Entropy-Based Moment Reduction for
Vlasov--Poisson--Dougherty Dynamics}
\author{Shiheng Zhang\thanks{Department of Applied Mathematics, University of
Washington, Seattle, WA 98195, USA (\texttt{shzhang3@uw.edu}).}
\and
Jiahao Zhang\thanks{Department of
Mathematics, Purdue University, 
West Lafayette, IN 47906, USA (\texttt{zhan2296@purdue.edu}).}
\and
Guang Lin\thanks{School of Mechanical Engineering and Department of
Mathematics, Purdue University, West Lafayette, IN 47906, USA
(\texttt{guanglin@purdue.edu}). Corresponding author.}}
\date{}

\maketitle

\begin{abstract}
We develop a data-adaptive entropy-based moment method for the
Vlasov--Poisson--Dougherty system with self-consistent temperature.
A positive exponential representation augments the physical moments
with data-adapted velocity functions.  The collision guarantees are independent of
basis construction: for every fixed basis satisfying the stated rank
condition, the reduced collision step is uniquely solvable, positive,
locally conservative and entropy-dissipative.  These guarantees
persist under a positive empirical cubature method (ECM), enabling
structure-preserving collision hyper-reduction.  We construct the basis by proper
orthogonal decomposition (POD), fit it with a neural network and
post-train it through the reduced dynamics.  In two-stream
relaxation, sustained heating and drive switching, one fixed
post-trained basis with thirteen variables per cell attains
space--time distribution errors below $0.25\%$ against same-grid
full-order solutions.  Post-training lowers this error relative to
the neural fit in all three tests and to POD in the two driven tests.
With full-grid MUSCL retained, matched drive-switching timings give
$\MatchedDynamicSpeedMinPrecise$--$\MatchedDynamicSpeedMaxPrecise$-fold
online speedups over the full-order solver using $128$--$\MatchedGridMax$
velocity nodes.
\end{abstract}



\section{Introduction}\label{sec:introduction}
The Vlasov--Poisson--Dougherty (VPD) system describes kinetic transport
and self-consistent electric fields coupled to nonlinear collisions.
The collision operator conserves density, momentum and kinetic energy,
and relaxes the distribution toward a local Maxwellian with its own
mean velocity and temperature.  Transport and forcing can continually
regenerate non-Gaussian velocity structure.  A useful reduced model
must capture this competition while retaining the collision invariants
and entropy dissipation.

Moment methods replace the velocity dependence by evolution equations
for finitely many velocity integrals.
Levermore's entropy-based hierarchy~\cite{Levermore1996}
uses constrained entropy minimization to obtain positive exponential
reconstructions and hyperbolic moment systems with an entropy.
The choice of moment functions determines which nonequilibrium
velocity structures these reconstructions can represent.

Data-driven methods address both the choice of moments and the cost
of closure evaluation.
Han et al.~\cite{Han2019} learn generalized moments together with
closures for their fluxes and collision terms.
Porteous et al.~\cite{PorteousLaiuHauck2023} approximate the
moment-to-entropy map by convex splines and neural networks;
Schotth\"ofer et al.~\cite{SchotthoferLaiuFrankHauck2025} develop
neural approximations of regularized entropy closures for radiation
transport.  These entropy approximations reduce the cost of closure
evaluation for prescribed moment functions.
Lara Benitez et al.~\cite{LaraBenitez2026} instead learn velocity
functions and coefficient maps for a discrete kinetic equilibrium.
Here, learned velocity functions describe the evolving nonequilibrium
distribution, and their coefficients are advanced by a numerical
collision solve.

For nonlinear collisions, reduction must also preserve the discrete
relations between the collision operator, its invariants and entropy.
Hakim et al.~\cite{Hakim2020Dougherty} construct conservative
discontinuous Galerkin schemes for the Dougherty operator.
Huang and Wang~\cite{HuangWang2026} formulate it as a projected
entropy gradient flow and use this structure in implicit particle
methods.  Related approaches include Landau discretizations with
positivity and exact H-theorems~\cite{HirvijokiBurbyKraus2018} and
entropy-dissipative scalar auxiliary variable schemes for kinetic
equations~\cite{ZhangShenHu2025}.
Micro--macro POD~\cite{ZhangHu2026MicroROM} and dynamical low-rank
methods~\cite{ZhangHu2026DLR} reduce VPFP systems with linear
Fokker--Planck collisions.  The VPD model considered here
has a temperature determined by the evolving distribution and
requires a reduced collision step that conserves all three physical
moments.

We construct such a step within a fixed exponential family that
retains the physical moment functions and augments them with
data-adapted velocity functions.  A projected velocity-grid mobility
and a strictly convex potential give a uniquely solvable, positive,
conservative and entropy-dissipative collision update for every
admissible fixed basis.  For a common input and frozen mobility,
the potential also characterizes the reduced update as a variational
best approximation to the full collision step.  These results
separate the design of the velocity basis from the collision
guarantees.

Snapshot reconstruction and dynamical prediction impose different
approximation objectives.  Otto et al.~\cite{OttoPadovanRowley2022}
optimize Petrov--Galerkin projections using trajectories of nonlinear
systems.  We initialize the velocity functions by proper orthogonal
decomposition (POD) of logarithmic snapshots, fit them with a neural
network and post-train through the implicit collision equations and
coupled transport.  The resulting basis stays fixed during prediction.
An exact KL decomposition separates representation error from
evolution within the family and helps interpret the effect of
post-training.

Reducing the state dimension does not by itself reduce the velocity
sums needed to assemble nonlinear operators.
Empirical cubature~\cite{HernandezCaicedoFerrer2017} selects integration
points and positive weights for reduced operators.
Conservative finite-volume reduction~\cite{CarlbergChoiSargsyan2018}
and entropy-stable hyper-reduction for conservation laws~\cite{Chan2020}
show how structural constraints can be retained during this stage.
We construct positive cubature rules for the Dougherty collision and
Newton matrices.  The sampled collision matrix annihilates the
collision invariants and remains positive semidefinite, so the same
variational proof gives the collision guarantees and entropy identity
on the original velocity grid.  Moment residuals and MUSCL transport
continue to use that grid.

\Needspace{3\baselineskip}
\paragraph{Contributions.}
\begin{itemize}
\setlength{\itemsep}{0pt}
\setlength{\parsep}{0pt}
\item We develop a low-dimensional, data-adaptive moment model that
captures non-Gaussian VPD dynamics with self-consistent temperature.
\item We prove that its implicit collision step is uniquely solvable,
positive, conservative and entropy-dissipative for every admissible
fixed basis.
\item We make the fixed velocity basis trainable through the reduced
dynamics and distinguish its representation and evolution errors
by an exact KL decomposition.
\item We construct positive empirical cubature that reduces collision
assembly while retaining the same invariants and discrete entropy
identity.
\end{itemize}

\section{The full-order model}\label{sec:full}
\subsection{The Vlasov--Poisson--Dougherty system}\label{sec:vpd}
The distribution $f(t,x,v)$ of a single kinetic species on periodic
$x\in[0,2\pi)$ and $v\in\R$ has the density, the current, the mean
velocity and the temperature
\begin{equation}\label{eq:moments}
 \rho=\int f\,dv,\qquad J=\rho u_f=\int vf\,dv,\qquad
 \rho T_f=\int(v-u_f)^2f\,dv .
\end{equation}
With a uniform neutralizing background of unit density, $f$ satisfies
\begin{align}
 \partial_tf+v\,\partial_xf+a\,\partial_vf&=Q(f),\qquad a=E+E_{\rm ext},
 \label{eq:vpd}\\
 \partial_xE&=\rho-1,\qquad \int_0^{2\pi}E\,dx=0 ,
 \label{eq:poisson}
\end{align}
where $E$ and $E_{\rm ext}(t,x)$ are the self-consistent and the
prescribed field, scaled by the signed charge-to-mass ratio so that $a$
is the acceleration, and $Q$ is the Dougherty operator with collision
frequency $\nu>0$~\cite{LenardBernstein1958,Dougherty1964},
\begin{equation}\label{eq:dougherty}
 Q(f)=\nu\,\partial_v\bigl((v-u_f)f+T_f\,\partial_vf\bigr).
\end{equation}
We take $\int\rho\,dx=2\pi$, so that~\eqref{eq:poisson} is solvable.
$Q$ conserves mass, momentum and energy, dissipates entropy and vanishes
exactly at the local Maxwellian
$M_{\rho,u_f,T_f}$, where
\begin{equation}\label{eq:maxwellian}
 M_{\rho,u,T}(v)=\frac{\rho}{\sqrt{2\pi T}}
 \exp\Bigl(-\frac{(v-u)^2}{2T}\Bigr),\qquad T>0 .
\end{equation}

Collisions therefore drop out of the equations for $\rho$, $J$ and the
kinetic energy.  Multiplying~\eqref{eq:vpd} by $1$, $v$ and $v^2/2$ and
integrating in $v$ gives
\begin{equation}\label{eq:moment-system}
 \begin{aligned}
 &\partial_t\rho+\partial_xJ=0,\qquad
 \partial_tJ+\partial_x\bigl(\rho u_f^2+\rho T_f\bigr)=a\rho,\\
 &\partial_t\bigl(\tfrac12\rho u_f^2+\tfrac12\rho T_f\bigr)
 +\tfrac12\partial_x\int v^3f\,dv=aJ,
 \end{aligned}
\end{equation}
and the energy flux contains the third central moment $S_3$:
\begin{equation}\label{eq:heat-flux}
 \int v^3f\,dv=\rho u_f^3+3\rho u_fT_f+S_3,\qquad
 S_3=\int(v-u_f)^3f\,dv .
\end{equation}
The third central moment $S_3$ is not determined by $\rho$, $u_f$ and
$T_f$.  By~\eqref{eq:vpd}
and~\eqref{eq:dougherty}, its evolution contains the source
$-3\rho T_f\,\partial_xT_f$ at a local Maxwellian and the collision term
$-3\nu S_3$: a temperature that varies in $x$ generates $S_3$, which
collisions damp only at the rate $3\nu$.

\subsection{Splitting and transport discretization}\label{sec:transport}
We discretize~\eqref{eq:vpd} by splitting it into
\begin{equation}\label{eq:split-subproblems}
 \begin{aligned}
 \partial_tf+v\,\partial_xf+a\,\partial_vf&=0
 &&\text{(transport)},\\
 \partial_tf&=Q(f)&&\text{(collisions)}.
 \end{aligned}
\end{equation}
A time step of length $\Delta t$ has the sequence
\[
 f^n\xrightarrow[\Delta t/2]{\text{collisions}}f^\star
 \xrightarrow[\Delta t]{\text{transport}}f^{\star\star}
 \xrightarrow[\Delta t/2]{\text{collisions}}f^{n+1}.
\]
We use $N_x$ spatial cells and $N_v\ge3$ velocity cells on $[-12,12]$,
\[
 \Delta x=\frac{2\pi}{N_x},\quad \Delta v=\frac{24}{N_v},\qquad
 x_i=(i-\tfrac12)\Delta x,\qquad v_j=-12+(j-\tfrac12)\Delta v ,
\]
with zero velocity flux at the boundary.  We write $\one$, $v$ and $v^2$
also for the vectors with entries $1$, $v_j$ and $v_j^2$, and $\exp$ and
$\log$ act componentwise.  For a positive nodal distribution
$f_{i,j}\approx f(t,x_i,v_j)$,
\begin{equation}\label{eq:density-shape}
 \rho_i=\Delta v\sum_jf_{i,j},\qquad
 q_{i,j}=\frac{\Delta v\,f_{i,j}}{\rho_i},\qquad \one^\top q_i=1 .
\end{equation}
In the transport subproblem of~\eqref{eq:split-subproblems}, the field
$E_i$ is computed from $\rho_i$ by solving~\eqref{eq:poisson} with the
fast Fourier transform.  With $a_i=E_i+E_{\rm ext}(t,x_i)$, the
semidiscrete transport equation is
\begin{equation}\label{eq:transport-rhs}
 \begin{aligned}
 \frac{df_{i,j}}{dt}=R(f,t)_{i,j}
 &=-\frac{\mathcal F^x_{i+1/2,j}-\mathcal F^x_{i-1/2,j}}{\Delta x}
 -\frac{\mathcal F^v_{i,j+1/2}-\mathcal F^v_{i,j-1/2}}{\Delta v},\\
 \mathcal F^v_{i,1/2}&=\mathcal F^v_{i,N_v+1/2}=0.
 \end{aligned}
\end{equation}
The MUSCL fluxes~\cite{VanLeer1979} use upwinding with limited linear
reconstructions and periodic boundaries in $i$.

\subsection{Structure-preserving collision step}\label{sec:mobility}
The collision subproblem of~\eqref{eq:split-subproblems} acts
independently in each spatial cell.
For a probability vector $q>0$, its temperature, its entropy and the
collision invariants on the grid are
\[
 T_q=(v^2)^\top q-(v^\top q)^2,\qquad H(q)=\sum_jq_j\log q_j,\qquad
 C=[\one,v,v^2/2]\in\R^{N_v\times3}.
\]
Using $\partial_vf=f\,\partial_v\log f$, we write~\eqref{eq:dougherty} as
\begin{equation}\label{eq:dougherty-flux}
 Q(f)=\nu T_f\,\partial_v
 \left[f\left(\partial_v\log f-\frac{u_f-v}{T_f}\right)\right].
\end{equation}
Following~\cite{HuangWang2026}, we discretize the gradient and its
weight at the cell interfaces
$v_{j+1/2}=-12+j\Delta v$:
\[
 (Dy)_{j+1/2}=\frac{y_{j+1}-y_j}{\Delta v},\qquad G=D[v,v^2/2],\qquad
 W_q=\diag\bigl(\Lambda(q_j,q_{j+1})\bigr),\qquad
 \Lambda(a,b)=\frac{b-a}{\log b-\log a},
\]
with $\Lambda(a,a)=a$ and $D\in\R^{(N_v-1)\times N_v}$.
The columns of $G$ sample $1$ and $v$ at the interfaces.  The
$W_q$-weighted projection onto the column space of $G$ is
\[
 P_q=G(G^\top W_qG)^{-1}G^\top W_q.
\]
For interface fluxes $\xi_{j+1/2}$, the velocity divergence is discretized as
\[
 (-D^\top\xi)_j
 =\frac{\xi_{j+1/2}-\xi_{j-1/2}}{\Delta v},\qquad
 \xi_{1/2}=\xi_{N_v+1/2}=0.
\]
For the collision step in~\eqref{eq:split-subproblems},
\eqref{eq:density-shape} and the flux discretization
of~\eqref{eq:dougherty-flux} give
\[
 \begin{aligned}
 f&=\frac{\rho}{\Delta v}q,\qquad
 \frac{d\rho}{dt}=0,\qquad D\log f=D\log q,\\
 \frac{\rho}{\Delta v}\frac{dq}{dt}
 &=\frac{df}{dt}
 =-\nu T_qD^\top\left[
 \frac{\rho}{\Delta v}W_q(I-P_q)D\log q\right].
 \end{aligned}
\]
Hence
\begin{align}
 \frac{dq}{dt}
 &=-\nu T_qD^\top W_q(I-P_q)D\log q,
 \label{eq:semidiscrete-collision}\\
 &=-\nu T_qK_q\log q,\qquad K_q=D^\top W_q(I-P_q)D.
 \label{eq:mobility}
\end{align}

\begin{proposition}\label{prop:discrete-collision}
Let $N_v\ge3$ and let $q$ be a probability vector with $q_j>0$.
The matrix $K_q$ in~\eqref{eq:mobility} satisfies
\begin{equation}\label{eq:invariants}
 K_q=K_q^\top\succeq0,\qquad K_qC=0,\qquad
 \ker K_q=\range(C).
\end{equation}
Along a positive solution of~\eqref{eq:semidiscrete-collision},
\[
 \frac{d}{dt}(C^\top q)=0,\qquad
 \frac{d}{dt}H(q)=-\nu T_q(\log q)^\top K_q\log q\le0.
\]
Thus collisions conserve mass, momentum, energy and $T_q$, and
dissipate entropy.  The positive equilibria are exactly the probability
vectors with $\log q\in\range(C)$; they include the normalized sampled
Maxwellians.
\end{proposition}

\begin{proof}
Since $W_q\succ0$ and $\rank G=2$, the definition of $P_q$ gives
\[
 P_q^2=P_q,\qquad P_qG=G,\qquad P_q^\top W_q=W_qP_q.
\]
Hence
\[
 K_q=\bigl[(I-P_q)D\bigr]^\top W_q\bigl[(I-P_q)D\bigr]\succeq0,
 \qquad K_qC=D^\top W_q(I-P_q)[0,G]=0.
\]
Moreover,
\[
 \ker K_q=\{y:Dy\in\range(G)\}
 =\spn\{\one,v,v^2/2\},
\]
because $\ker D=\spn\{\one\}$ and $G=D[v,v^2/2]$.
Multiplying~\eqref{eq:semidiscrete-collision} by $C^\top$ and
$(\one+\log q)^\top$ gives conservation and entropy dissipation.
Finally, $\nu T_q>0$, so the collision right-hand side vanishes
exactly when $\log q\in\ker K_q=\range(C)$.
\end{proof}

\begin{remark}
For smooth positive $f$ with vanishing velocity boundary terms,
integration by parts gives
\[
 \int_{\R} f(v)\begin{pmatrix}1\\v\end{pmatrix}
 \left[\partial_v\log f-\frac{u_f-v}{T_f}\right]\,dv=0.
\]
The matrix $P_q$ enforces the corresponding interface orthogonality.
On a truncated grid, this correspondence is affected by tail and
discretization errors, examined in Section~\ref{sec:refinement} and
the analytic benchmark in Appendix~\ref{app:implementation}.
\end{remark}

For a collision substep of length $h>0$, the implicit update is
\begin{equation}\label{eq:full}
 \frac{q^{n+1}-q^n}{h}=-\nu T_{q^n}K_{q^n}\log q^{n+1},
\end{equation}
where $q^n$ and $q^{n+1}$ are the cell probability vectors before and
after this substep.  By Proposition~\ref{prop:collision},
\eqref{eq:full} has a unique solution $q^{n+1}>0$, with
$C^\top q^{n+1}=C^\top q^n$ and $H(q^{n+1})\le H(q^n)$.
Applied cellwise, this update defines
\[
 \mathcal C_h(f)_{i,j}=\frac{\rho_i}{\Delta v}\,q_j^{n+1},
 \qquad q^n=q_i.
\]

\subsection{Coupled time step}\label{sec:full-step}
Combining the collision update $\mathcal C_h$ with the transport
right-hand side $R$ gives one step from $f^n$ at $t_n=n\Delta t$ to
$f^{n+1}$:
\begin{equation}\label{eq:splitting}
 \begin{aligned}
 f^\star&=\mathcal C_{\Delta t/2}(f^n),\\
 g&=f^\star+\Delta t\,R(f^\star,t_n),\\
 f^{\star\star}&=\tfrac12f^\star+\tfrac12\bigl[g+\Delta t\,R(g,t_n+\Delta t)\bigr],\\
 f^{n+1}&=\mathcal C_{\Delta t/2}(f^{\star\star}),
 \end{aligned}
\end{equation}
with $f^0_{i,j}=f_0(x_i,v_j)$ for the initial distribution $f_0>0$.
The two middle stages form the two-stage strong-stability-preserving
(SSP) Runge--Kutta transport step, with the field computed from the
state of each stage.  We require
$\Delta t\,\bigl(\max_j|v_j|/\Delta x+\max_i|a_i|/\Delta v\bigr)\le0.48$ at
both transport stages.  Under this bound, every stage
of~\eqref{eq:splitting} stays positive (Appendix~\ref{app:implementation}).
The first-order collision update makes the coupled method first order
in $\Delta t$ (Section~\ref{sec:refinement}).

\section{Entropy-based moment model}\label{sec:reduced}
We construct a moment approximation $\mathcal C_h^{\rm red}$ of the
collision map $\mathcal C_h$ in~\eqref{eq:splitting}.

\subsection{Exponential family and moment matching}\label{sec:family}\label{sec:projection}
The logarithm of a Maxwellian~\eqref{eq:maxwellian} is quadratic,
\[
 \log M_{\rho,u,T}(v)=\mathrm{const}+\frac uT\,v-\frac{v^2}{2T}.
\]
We retain $1,v,v^2$ and add $d$ fixed functions $\psi_1,\ldots,\psi_d$ to
represent non-Gaussian distributions.  We assume that
$1,v,v^2,\psi_1,\ldots,\psi_d$ are linearly independent on the velocity
grid.  Let $B\in\R^{N_v\times r}$, $r=d+2$, be the fixed basis matrix
constructed from these functions in Section~\ref{sec:basis}, with row
$B_j$ at $v_j$.  The constant function is omitted from $B$, since its
exponential factor cancels in the normalization below.
For coefficients $z_i\in\R^r$, define
\begin{equation}\label{eq:family}
 \begin{aligned}
 \ell_i(v_j)&=B_jz_i=\sum_{k=1}^{r}B_{jk}z_{i,k},\\
 p_j(z_i)&=\frac{e^{\ell_i(v_j)}}{Z(z_i)},\qquad
 Z(z_i)=\sum_{s=1}^{N_v}e^{\ell_i(v_s)},\\
 f_{i,j}&=\frac{\rho_i}{\Delta v}\,p_j(z_i).
 \end{aligned}
\end{equation}
At each time level, the distribution is represented by $(\rho_i,z_i)$,
with the basis fixed.

Given a cell probability vector $q_i>0$, we match its moments by
\begin{align}
 m(z_i)&=B^\top p(z_i),\label{eq:moment-map}\\
 m(z_i)&=B^\top q_i
 \quad\Longleftrightarrow\quad
 \sum_jB_{jk}p_j(z_i)=\sum_jB_{jk}q_{i,j},
 \quad k=1,\ldots,r.\label{eq:moment-matching}
\end{align}
Since $1,v,v^2$ are retained, this gives $C^\top p(z_i)=C^\top q_i$:
with $\rho_i$ unchanged, the reconstruction preserves the density, mean
velocity and temperature.  It also retains the moments of
$\psi_1,\ldots,\psi_d$.

\subsection{Reduced collision step and time stepping}\label{sec:reduced-collision}\label{sec:newton}
For the collision step~\eqref{eq:full} with $h>0$, set $q^{n+1}=p(z)$.
Since $K_{q^n}\one=0$,
\[
 K_{q^n}\log p(z)
 =K_{q^n}\bigl[Bz-\log Z(z)\one\bigr]=K_{q^n}Bz.
\]
We require the collision residual to have zero moments:
\[
 B^\top\bigl[p(z)-q^n+h\nu T_{q^n}K_{q^n}\log p(z)\bigr]=0.
\]
With $L_q=B^\top K_qB\succeq0$, this is the $r$-dimensional equation
\begin{equation}\label{eq:reduced}
 F(z):=m(z)-B^\top q^n+h\nu T_{q^n}L_{q^n}z=0.
\end{equation}
Differentiating~\eqref{eq:moment-map} gives
\begin{equation}\label{eq:moment-jacobian}
 \nabla m(z)=B^\top\bigl[\diag p(z)-p(z)p(z)^\top\bigr]B.
\end{equation}
We solve~\eqref{eq:reduced} by damped Newton iteration.  The Newton
correction $\delta z\in\R^r$ and the coefficient update are
\[
 \bigl[\nabla m(z)+h\nu T_{q^n}L_{q^n}\bigr]\delta z=-F(z),
 \qquad z\leftarrow z+\alpha\,\delta z,
\]
with the step length $0<\alpha\le1$ chosen by line search.
The collision update in each cell is
\[
 q^{n+1}=p(z),\qquad
 \mathcal C_h^{\rm red}(f)_{i,j}
 =\frac{\rho_i}{\Delta v}\,q_j^{n+1}.
\]
The coupled scheme is~\eqref{eq:splitting} with initial data
reconstructed by~\eqref{eq:moment-matching} and $\mathcal C_h$
replaced by $\mathcal C_h^{\rm red}$.

\subsection{Solvability, conservation and entropy}\label{sec:theorem}
\begin{proposition}\label{prop:collision}
Let $h>0$ and let $q^n>0$ satisfy $\one^\top q^n=1$.
Under the basis assumptions of Section~\ref{sec:family},
\eqref{eq:reduced} has a unique solution $z\in\R^r$.
The probability vector $q^{n+1}=p(z)>0$ satisfies
\begin{equation}\label{eq:reduced-entropy}
 \begin{gathered}
 C^\top q^{n+1}=C^\top q^n,\\
 H(q^{n+1})-H(q^n)+\KL(q^n\Vert q^{n+1})
 +h\nu T_{q^n}\,z^\top L_{q^n}z=0.
 \end{gathered}
\end{equation}
Here $\KL(a\Vert b)=\sum_ja_j\log(a_j/b_j)$, so
$H(q^{n+1})\le H(q^n)$.
The full update~\eqref{eq:full} also has a unique positive
solution, conserves $C^\top q$ and satisfies $H(q^{n+1})\le H(q^n)$.
\end{proposition}

\begin{proof}
Let
\begin{equation}\label{eq:rank}
 S_B=\spn\{\one,B_{:,1},\ldots,B_{:,r}\}
 =\spn\{1,v,v^2,\psi_1,\ldots,\psi_d\}.
\end{equation}
For $\tau\ge0$, a fixed matrix $K=K^\top\succeq0$ with $KC=0$ and a
subspace $S$ containing $\range(C)$, consider
\begin{equation}\label{eq:potential}
 \begin{aligned}
 \mathcal V(y)&=\sum_j\bigl(e^{y_j}-q^n_jy_j\bigr)
                  +\frac\tau2\,y^\top Ky,\\
 \nabla\mathcal V(y)&=e^y-q^n+\tau Ky,\qquad
 \nabla^2\mathcal V(y)=\diag(e^y)+\tau K\succ0.
 \end{aligned}
\end{equation}
Each term $e^{y_j}-q^n_jy_j$ is bounded below by
$q^n_j(1-\log q^n_j)$ and tends to $+\infty$ as $|y_j|\to\infty$.
Thus $\mathcal V$ is coercive and strictly convex, and has a unique
minimizer $y$ on $S$.  With $q^{n+1}=e^y>0$, its stationarity condition is
\begin{equation}\label{eq:galerkin}
 \varphi^\top\bigl(q^{n+1}-q^n+\tau K\log q^{n+1}\bigr)=0
 \qquad\text{for all }\varphi\in S.
\end{equation}
Testing with the columns of $C$, using $C^\top K=0$, gives
\[
 C^\top q^{n+1}=C^\top q^n,\qquad \one^\top q^{n+1}=1.
\]
Testing with $\varphi=\log q^{n+1}\in S$ and using
\[
 H(q^{n+1})-H(q^n)+\KL(q^n\Vert q^{n+1})
 =(q^{n+1}-q^n)^\top\log q^{n+1}
\]
gives
\begin{equation}\label{eq:entropy-identity}
 H(q^{n+1})-H(q^n)+\KL(q^n\Vert q^{n+1})
 +\tau\,(\log q^{n+1})^\top K\log q^{n+1}=0.
\end{equation}

Set $\tau=h\nu T_{q^n}$ and $K=K_{q^n}$, whose matrix properties follow
from Proposition~\ref{prop:discrete-collision}.
For $h>0$ and $S=\R^{N_v}$, \eqref{eq:galerkin} is the full
update~\eqref{eq:full} multiplied by $h$.
For $S=S_B$, linear independence gives a unique representation
$y=\lambda\one+Bz$.  Testing in the directions $\one$ and $B_{:,k}$ gives
\[
 \begin{aligned}
 e^\lambda Z(z)&=1,\qquad q^{n+1}=p(z),\\
 B^\top\nabla\mathcal V(y)
 &=m(z)-B^\top q^n+\tau L_{q^n}z=F(z)=0.
 \end{aligned}
\]
Thus stationarity on $S_B$ is equivalent to~\eqref{eq:reduced}, giving
its unique solution.  The Newton matrix is positive definite since
$L_{q^n}\succeq0$ and~\eqref{eq:moment-jacobian} gives
\[
 \xi^\top\nabla m(z)\,\xi
 =\sum_jp_j(z)\bigl(B_j\xi-m(z)^\top\xi\bigr)^2>0
 \qquad(\xi\ne0).
\]
Equation~\eqref{eq:entropy-identity} gives~\eqref{eq:reduced-entropy},
using $K_{q^n}\one=0$.
Minimizing $\mathcal V$ in~\eqref{eq:potential} on $S_B$ with $\tau=0$
gives existence and uniqueness for~\eqref{eq:moment-matching}.
\end{proof}

The reconstruction~\eqref{eq:moment-matching} minimizes entropy:
\[
 p(z)=\arg\min\bigl\{H(p'):\ p'\ge0,\ \one^\top p'=1,\ B^\top p'=B^\top q^n\bigr\}.
\]
For every feasible $p'$, moment matching and $\log p(z)\in S_B$ give
$H(p')-H(p(z))=\KL(p'\Vert p(z))\ge0$, with $0\log0=0$.
The collision guarantees apply to every fixed basis satisfying the
assumptions of Section~\ref{sec:family}.
Since $\range(C)\subset S_B$, the family also contains every positive
equilibrium of~\eqref{eq:semidiscrete-collision}.

\subsection{Approximation properties}\label{sec:approximation}
The variational formulation~\eqref{eq:potential} also characterizes
which distribution the reduced collision step selects.

\begin{proposition}\label{prop:collision-approximation}
Under the assumptions of Proposition~\ref{prop:collision}, let
$q_F^{n+1}$ and $p(z_R)$ solve~\eqref{eq:full}
and~\eqref{eq:reduced}, respectively, with the same $q^n$ and $h>0$.
Then
\begin{equation}\label{eq:collision-best-approximation}
 z_R=\underset{z\in\R^r}{\arg\min}\left\{
 \KL\bigl(q_F^{n+1}\Vert p(z)\bigr)
 +\frac{h\nu T_{q^n}}2
 \bigl[\log p(z)-\log q_F^{n+1}\bigr]^\top K_{q^n}
 \bigl[\log p(z)-\log q_F^{n+1}\bigr]\right\}.
\end{equation}
\end{proposition}

\begin{proof}
Use $\tau=h\nu T_{q^n}$ and $K=K_{q^n}$
in~\eqref{eq:potential}.  Since
$\nabla\mathcal V(\log q_F^{n+1})=0$, direct expansion gives
\[
 \begin{aligned}
 &\mathcal V(\log p(z))-\mathcal V(\log q_F^{n+1})\\
 &\quad=\KL\bigl(q_F^{n+1}\Vert p(z)\bigr)\\
 &\qquad\quad+\frac\tau2
 \bigl[\log p(z)-\log q_F^{n+1}\bigr]^\top K
 \bigl[\log p(z)-\log q_F^{n+1}\bigr].
 \end{aligned}
\]
This is the objective in~\eqref{eq:collision-best-approximation},
uniquely minimized at $z_R$ by Proposition~\ref{prop:collision}.
\end{proof}

The KL error of the normalized cell distributions can be separated
by the exponential-family Pythagorean identity~\cite{Csiszar1975},
using moment matching~\eqref{eq:moment-matching} at each time.
\begin{proposition}\label{prop:kl-decomposition}
At a fixed time and cell $i$, let $q_i^F>0$ be the full-model
probability vector, let $p(z_i^*)$ match its moments by
\eqref{eq:moment-matching}, and let $p(z_i^R)$ be any reduced-model
probability vector in the same family.  Then
\begin{equation}\label{eq:kl-error-decomposition}
 \KL\bigl(q_i^F\Vert p(z_i^R)\bigr)
 =\underbrace{\KL\bigl(q_i^F\Vert p(z_i^*)\bigr)}_{\text{representation error}}
 +\underbrace{\KL\bigl(p(z_i^*)\Vert p(z_i^R)\bigr)}_{\text{error within the family}}.
\end{equation}
The first term is the minimum KL error attainable in the fixed family.
\end{proposition}

\begin{proof}
The difference between the left side and the two terms on the right is
\[
 \begin{aligned}
 &(q_i^F-p(z_i^*))^\top
       \bigl[\log p(z_i^*)-\log p(z_i^R)\bigr]\\
 &\quad=(q_i^F-p(z_i^*))^\top
 \left[B(z_i^*-z_i^R)
       -\log\frac{Z(z_i^*)}{Z(z_i^R)}\one\right]=0,
 \end{aligned}
\]
by moment matching and unit mass.  Both terms on the right of
\eqref{eq:kl-error-decomposition} are nonnegative, and the second is
zero exactly when $p(z_i^R)=p(z_i^*)$.
\end{proof}

The error within the family in~\eqref{eq:kl-error-decomposition}
measures how closely the evolved state follows the best reconstruction;
it can include closure, cubature and discretization effects.
This normalized diagnostic excludes density errors.
For a collision input $q^n$, let $p(z^*)$ match its moments.
Then~\eqref{eq:kl-error-decomposition} splits the entropy decrease
in~\eqref{eq:reduced-entropy} as
\[
 \begin{aligned}
 H(q^n)-H(q^{n+1})
 &=\KL(q^n\Vert p(z^*))+\KL(p(z^*)\Vert q^{n+1})\\
 &\quad+h\nu T_{q^n}z^\top L_{q^n}z.
 \end{aligned}
\]
The first term is the decrease due to moment reconstruction; the
remaining terms arise from the implicit collision update.
Transport entropy changes are measured separately in
Section~\ref{sec:structure}.

\section{Structure-preserving empirical cubature}\label{sec:quadrature}\label{sec:rules}
Assembling the matrices in~\eqref{eq:reduced}
and~\eqref{eq:moment-jacobian} by exact sums costs $O(N_vr^2)$ per cell.
We use a greedy empirical cubature method (ECM)~\cite{Chan2020,HernandezCaicedoFerrer2017},
with nonnegative least-squares (NNLS) weights, to construct
\[
 L_{q^n}\longrightarrow\widetilde L_{q^n},\qquad
 \nabla m(z)\longrightarrow J_{\mathcal N}(z).
\]
The rules are fitted offline for each fixed basis and velocity grid;
Appendix~\ref{app:implementation} gives the fitting details.

\subsection{Quadrature for the collision matrix}\label{sec:edge-quadrature}
Using $K_q$ from~\eqref{eq:mobility},
\[
 L_q=B^\top K_qB
 =\bigl[(I-P_q)DB\bigr]^\top W_q\bigl[(I-P_q)DB\bigr].
\]
Let $\Gamma\subset\{1,\ldots,N_v-1\}$ contain $N_\Gamma$ selected
interfaces with weights $\alpha_e>0$.  Let $D_\Gamma$ and $G_\Gamma$
contain the corresponding rows of $D$ and $G$.
For $q>0$ and $\rank G_\Gamma=2$, define
\begin{equation}\label{eq:sparse-edge-metric}
 \begin{aligned}
 W_\Gamma(q)&=\diag\bigl(\alpha_e\Lambda(q_e,q_{e+1})\bigr)_{e\in\Gamma},\\
 P_\Gamma&=G_\Gamma(G_\Gamma^\top W_\Gamma G_\Gamma)^{-1}
               G_\Gamma^\top W_\Gamma,\\
 \widetilde L_q&=\bigl[(I-P_\Gamma)D_\Gamma B\bigr]^\top
                  W_\Gamma\bigl[(I-P_\Gamma)D_\Gamma B\bigr].
 \end{aligned}
\end{equation}
The collision equation~\eqref{eq:reduced} becomes
\begin{equation}\label{eq:sparse-collision-residual}
 \widetilde F(z):=m(z)-B^\top q^n+h\nu T_{q^n}\widetilde L_{q^n}z=0,
\end{equation}
where $p(z)$, $m(z)$, $B^\top q^n$ and $T_{q^n}$ are evaluated on the
full velocity grid.

\begin{proposition}\label{prop:sparse-collision}
Under the assumptions of Proposition~\ref{prop:collision}, let $\alpha_e>0$
for $e\in\Gamma$ and $\rank G_\Gamma=2$.
Then~\eqref{eq:sparse-collision-residual} has a unique solution
$z\in\R^r$.  The probability vector $q^{n+1}=p(z)>0$ satisfies
$C^\top q^{n+1}=C^\top q^n$ on the full velocity grid and
\begin{equation}\label{eq:sparse-entropy}
 H(q^{n+1})-H(q^n)+\KL(q^n\Vert q^{n+1})
 +h\nu T_{q^n}\,z^\top\widetilde L_{q^n}z=0.
\end{equation}
In particular, $H(q^{n+1})\le H(q^n)$.
\end{proposition}

\begin{proof}
Set
\[
 \begin{aligned}
 \widetilde K_q
 &=\bigl[(I-P_\Gamma)D_\Gamma\bigr]^\top
       W_\Gamma\bigl[(I-P_\Gamma)D_\Gamma\bigr]
   =\widetilde K_q^\top\succeq0,\\
 \widetilde L_q&=B^\top\widetilde K_qB,\\
 D_\Gamma C&=[0,G_\Gamma],\quad P_\Gamma G_\Gamma=G_\Gamma
 \quad\Longrightarrow\quad \widetilde K_qC=0.
 \end{aligned}
\]
The variational proof of Proposition~\ref{prop:collision} applies with
$K_{q^n}$ replaced by $\widetilde K_{q^n}$ and
gives~\eqref{eq:sparse-entropy}.
\end{proof}

An additional rank condition preserving the reduced equilibrium
nullspace is checked in Appendix~\ref{app:implementation}.
Proposition~\ref{prop:collision-approximation} also applies when both
the full and reduced collision steps use $\widetilde K_{q^n}$ in place
of $K_{q^n}$, with the same input $q^n$.

\subsection{Quadrature for the Newton matrix and online cost}\label{sec:node-quadrature}
For a set $\mathcal N$ of $N_{\mathcal N}$ selected velocity nodes with
weights $\alpha_j>0$, define
\[
 \pi_j(z)=\frac{\alpha_je^{B_jz}}
                   {\sum_{s\in\mathcal N}\alpha_se^{B_sz}},
 \qquad j\in\mathcal N.
\]
Approximating the covariance in~\eqref{eq:moment-jacobian} gives
\begin{equation}\label{eq:sparse-covariance}
 J_{\mathcal N}(z)=\sum_{j\in\mathcal N}\pi_j(z)B_j^\top B_j
 -\Bigl(\sum_{j\in\mathcal N}\pi_j(z)B_j\Bigr)^\top
  \Bigl(\sum_{j\in\mathcal N}\pi_j(z)B_j\Bigr).
\end{equation}
If $1,v,v^2,\psi_1,\ldots,\psi_d$ are linearly independent on
$\mathcal N$, then $J_{\mathcal N}(z)\succ0$, by the covariance argument
in Proposition~\ref{prop:collision}.
Replacing $\nabla m(z)$ in Section~\ref{sec:newton} gives a
safeguarded quasi-Newton iteration,
\begin{equation}\label{eq:sparse-newton}
 \bigl[J_{\mathcal N}(z)+h\nu T_{q^n}\widetilde L_{q^n}\bigr]\delta z
 =-\widetilde F(z).
\end{equation}
The line search and stopping tests evaluate the
residual~\eqref{eq:sparse-collision-residual} on the full velocity grid;
Appendix~\ref{app:implementation} specifies the safeguards.

For a fixed basis, the per-cell costs of the weighted sums are
\begin{center}
\small
\begin{tabular}{lcc}
\toprule
Quantity & Exact sums & ECM\\
\midrule
Collision matrix & $O(N_vr^2)$ & $O(N_\Gamma r^2)$\\
Newton covariance & $O(N_vr^2)$ & $O(N_{\mathcal N}r^2)$\\
Full-grid moments & $O(N_vr)$ & $O(N_vr)$\\
\bottomrule
\end{tabular}
\end{center}
The $r\times r$ linear solve costs $O(r^3)$ in both cases.
Normalization and MUSCL transport are evaluated on the full velocity grid.

\section{Basis construction}\label{sec:basis}
From training probability vectors $q^{(1)},\ldots,q^{(N_s)}$, we
construct the $d$ non-Gaussian functions in three stages:
\[
 \{\log q^{(n)}\}_{n=1}^{N_s}
 \;\xrightarrow{\ \text{POD}\ }\;
 \Psi_{\rm POD}
 \;\xrightarrow{\ \text{network fit}\ }\;
 \Psi_{\vartheta_0}
 \;\xrightarrow{\ \text{post-training}\ }\;
 \Psi_{\vartheta_*}.
\]
Only the functions $\psi_k$ change; $1,v,v^2$ remain analytic.
Each stage produces a fixed basis to which the guarantees of
Sections~\ref{sec:reduced} and~\ref{sec:quadrature} apply.

A weighted QR converts the samples $\Psi_{jk}=\psi_k(v_j)$ and
$1,v,v^2$ into $B$ without changing $S_B$ in~\eqref{eq:rank}.
Appendix~\ref{app:implementation} gives the coordinates.

\subsection{Training data and POD}\label{sec:pod}
The four training and two validation parameter sets are listed in
Table~\ref{tab:parameters}; Appendix~\ref{app:training} gives the grids
and sampling intervals.  For fixed weights $\omega_j>0$, set
\[
 \Omega=\diag(\omega_j),\qquad
 \Pi_C=C(C^\top\Omega C)^{-1}C^\top\Omega.
\]
The POD functions are
\begin{equation}\label{eq:pod}
 \begin{aligned}
 Y&=(I-\Pi_C)\bigl[\log q^{(1)},\ldots,\log q^{(N_s)}\bigr],\\
 \Omega^{1/2}Y&=U\Sigma V^\top,\\
 \Psi_{\rm POD}&=\Omega^{-1/2}U_{:,1:d},
 \end{aligned}
\end{equation}
where the singular values decrease along $\Sigma$.  Removing the
quadratic part represented by $1,v,v^2$ makes POD act on the
non-Gaussian part of $\log q$.  The truncated SVD gives its best
rank-$d$ approximation in the $\Omega$-weighted norm.  On other grids,
the POD columns are evaluated by natural cubic splines.

\subsection{Network approximation of the POD basis}\label{sec:neural-pod}
We represent the POD functions by a differentiable neural network
$\mathcal N_\vartheta:\R\to\R^d$, with
\[
 \psi_k(v;\vartheta)=\gamma_k(\mathcal N_\vartheta(v))_k.
\]
We fit the POD functions on the training grid by
\begin{equation}\label{eq:neural-pod}
 \min_\vartheta\ \frac1d\sum_{k=1}^{d}\sum_{j=1}^{N_v}
 w_j\,
 \frac{\bigl(\psi_k(v_j;\vartheta)-(\Psi_{\rm POD})_{jk}\bigr)^2}{\gamma_k^2}.
\end{equation}
The fixed output scales $\gamma_k>0$ and fitting weights $w_j>0$ are
specified in Appendix~\ref{app:training}.
Validation selects $\vartheta_0$ by the mean
$\KL(q\Vert p(z))$ after moment matching~\eqref{eq:moment-matching}.
The network is evaluated on other velocity grids without refitting.
Online, its parameters stay fixed and only $(\rho_i,z_i)$ evolve.

\subsection{Post-training through the dynamics}\label{sec:dynamics-training}
POD fits logarithmic snapshots, whereas the coupled method predicts a
trajectory.  We therefore start from $\vartheta_0$ and adjust the basis
for the prediction error.  For each $\vartheta$, we recompute the
weighted QR defining $B_\vartheta$ (Appendix~\ref{app:implementation}).

For each training interval $w\in\mathcal W$, use the saved full-order
solutions $f_w^n$ at $t_w+n\Delta t$, $n\in\{0,16,32\}$.
The reduced prediction starts from moment matching and follows the
coupled scheme with collision quadrature:
\[
 f_w^0\ \xrightarrow{\ \eqref{eq:moment-matching}\ }\
 \widehat f_{\vartheta,w}^{\,0}
 \ \xrightarrow{\ \eqref{eq:splitting},\,\eqref{eq:sparse-collision-residual}\ }\
 \widehat f_{\vartheta,w}^{\,1}
 \ \longrightarrow\ \cdots\ \longrightarrow\
 \widehat f_{\vartheta,w}^{\,32}.
\]
The node and interface rules fitted for the network approximation stay
fixed during post-training.  With the Frobenius norm $\|\cdot\|_F$, we minimize
\begin{equation}\label{eq:dynamic-loss}
 \begin{split}
 \mathcal L_{\rm dyn}(\vartheta)
 =\frac{1}{2|\mathcal W|}
 \sum_{w\in\mathcal W}\sum_{n\in\{16,32\}}
 \bigg[&\frac{\|\widehat f_{\vartheta,w}^{\,n}-f_w^n\|_F^2}
                  {\|f_w^n\|_F^2}
 +0.1\frac{\|\widehat J_{\vartheta,w}^{\,n}-J_w^n\|_F^2}{N_x\overline{J^2}}\\
 &+0.1\frac{\|\widehat S_{3,\vartheta,w}^{\,n}-S_{3,w}^n\|_F^2}
                  {N_x\overline{S_3^2}}\bigg],
 \end{split}
\end{equation}
where $J$ and $S_3$ are evaluated cellwise from~\eqref{eq:moments}
and~\eqref{eq:heat-flux}, and $\overline{J^2}$ and
$\overline{S_3^2}$ are their mean squares
over the training trajectories.  Figure~\ref{fig:basis-training} uses
\eqref{eq:kl-error-decomposition} to assess separately how post-training
changes the best reconstruction and the evolved trajectory.

Differentiating each converged root
$\widetilde F_\vartheta(z;q^n)=0$ of~\eqref{eq:sparse-collision-residual} gives
\[
 \partial_z\widetilde F_\vartheta\,\mathrm dz
 =-\partial_\vartheta\widetilde F_\vartheta\,\mathrm d\vartheta
  -\partial_{q^n}\widetilde F_\vartheta\,\mathrm dq^n.
\]
This derivative uses the exact residual Jacobian; the sampled matrix
in~\eqref{eq:sparse-newton} is used only in the forward iteration.
Appendix~\ref{app:training} gives the adjoint evaluation and training parameters.

After training, we freeze $\vartheta_*$, refit the quadrature rules and
check the basis and quadrature rank conditions on each grid.
The test problems of Section~\ref{sec:numerics} enter neither training
nor selection of $\vartheta_*$.

\section{Numerical experiments}\label{sec:numerics}
\subsection{Test problems and error measures}\label{sec:setup}
The training, validation and two-stream test trajectories start from
mixtures of two Gaussians with spatially varying weight and shift,
\begin{equation}\label{eq:initial}
 f_0(x,v)=(1+A_\rho\cos x)\bigl[\eta(x)\,M_{1,c+s(x),T_+}(v)
 +(1-\eta(x))\,M_{1,-c+s(x),T_-}(v)\bigr],
\end{equation}
where $\eta(x)=\eta_0+A_\eta\sin x$, $s(x)=A_s\sin x$, and each Maxwellian
is renormalized on the grid, so that $\rho_i=1+A_\rho\cos x_i$.  The
external field is
\begin{equation}\label{eq:external-field}
 E_{\rm ext}(t,x)=A_{\rm ext}\bigl[\sin(x-\omega t)
 +\tfrac12\sin(2x-1.7t+0.4)\bigr],
\end{equation}
with the parameters of Table~\ref{tab:parameters}.  The three two-stream
tests are unforced relaxation up to $t=4$, sustained
heating up to $t=8$ ($t=4$ in Section~\ref{sec:refinement}), and drive
switching up to $t=12$.  Drive switching reverses and then removes the
heating drive by replacing $A_{\rm ext}=0.7$ with
\begin{equation}\label{eq:switch-drive}
 A_{\rm ext}(t)=\begin{cases}
 0.7,&0\le t<4,\\
 -0.7,&4\le t<8,\\
 0,&8\le t\le12.
 \end{cases}
\end{equation}
Up to $t=4$, drive switching coincides with heating, and in step $n$ both
transport stages of~\eqref{eq:splitting} use $A_{\rm ext}(t_n)$.

\begin{table}[!b]
\centering\small
\setlength{\tabcolsep}{4pt}
\caption{Parameters of the initial data~\eqref{eq:initial}, the collision
frequency $\nu$ and the external field~\eqref{eq:external-field} for the
training and validation trajectories of Section~\ref{sec:pod} and the
test problems.}
\label{tab:parameters}
\begin{tabular}{lrrrrrrrrrr}
\toprule
Trajectory & $A_\rho$ & $c$ & $T_+$ & $T_-$ & $\eta_0$ & $A_\eta$ & $A_s$ & $\nu$ & $A_{\rm ext}$ & $\omega$\\
\midrule
Training A & 0.15 & 1.6 & 0.4 & 0.7 & 0.6 & 0.07 & 0.2 & 0.3 & 0.35 & 0.9 \\
Training B & 0.12 & 2 & 0.6 & 0.45 & 0.45 & 0.1 & 0.1 & 0.5 & 0.8 & 1.25 \\
Training C & 0.2 & 1 & 0.45 & 0.8 & 0.5 & 0.04 & 0.2 & 0.15 & 0.55 & 1.4 \\
Training D & 0.16 & 1.7 & 0.25 & 0.5 & 0.5 & 0.12 & 0.22 & 0.4 & 0.45 & 1.3 \\
Validation A & 0.13 & 1.5 & 0.55 & 0.4 & 0.52 & 0.05 & 0.13 & 0.25 & 0.6 & 1 \\
Validation B & 0.14 & 1.4 & 0.3 & 0.45 & 0.48 & 0.07 & 0.18 & 0.3 & 0.5 & 1.2 \\
\midrule
Two-stream relaxation & 0.18 & 1.8 & 0.35 & 0.6 & 0.55 & 0.08 & 0.15 & 0.35 & 0 & -- \\
Sustained heating & 0.1 & 1.2 & 0.5 & 0.5 & 0.5 & 0.06 & 0.1 & 0.2 & 0.7 & 1.1 \\
\bottomrule
\end{tabular}
\end{table}

Unless stated otherwise, $N_x=64$, $N_v=128$, $\Delta t=0.00125$ and
$d=10$, so that $r=12$ and each cell has thirteen state variables.
The reduced-order model uses the post-trained basis and positive ECM,
except in the rank comparison of Table~\ref{tab:rank} and the
refinement study of Section~\ref{sec:refinement}.
Outside Section~\ref{sec:refinement}, each reduced computation is
compared with the full-order model on the same grids and with the same time
step, so the reported errors are errors of the reduction, including that
of the initial reconstruction.

With $J$ and $S_3$ as in Section~\ref{sec:dynamics-training}, we write
$T_f=(T_i)_i$ for the temperatures $T_i=T_{q_i}$ of the cells, and use
$\|g\|^2=\Delta x\Delta v\sum_{i,j}g_{i,j}^2$ for $f$ and
$\|g\|^2=\Delta x\sum_ig_i^2$ for functions of $x$.  For a time interval $\mathcal I$, the space--time relative error
in $f$ is
\begin{equation}\label{eq:error}
 e_{\mathcal I}(f)=\biggl(\frac{\int_{\mathcal I}\|f_{\rm red}(t)-f_{\rm full}(t)\|^2\,dt}
 {\int_{\mathcal I}\|f_{\rm full}(t)\|^2\,dt}\biggr)^{1/2},
\end{equation}
where the time integrals use the trapezoidal rule on snapshots spaced
by $0.04$, and $e_{\mathcal I}(J)$, $e_{\mathcal I}(T_f)$ and
$e_{\mathcal I}(S_3)$ are defined in the same way.  Without a stated
interval, $\mathcal I$ is the whole run.  The instantaneous relative error is
$\|f_{\rm red}(t)-f_{\rm full}(t)\|/\|f_{\rm full}(t)\|$, and likewise for $J$
and $T_f$.  The mass-weighted mean temperature is
$\overline T_f=\sum_i\rho_iT_i/\sum_i\rho_i$, and the relative deviation from the local Maxwellian is
\begin{equation}\label{eq:nongaussian}
 d_G(t)=\frac{\|f-f_G\|}{\|f\|},\qquad
 (f_G)_{i,j}=\frac{\rho_i}{\Delta v}\,
 \frac{M_{1,u_i,T_i}(v_j)}{\one^\top M_{1,u_i,T_i}(v)} .
\end{equation}

\subsection{Non-Gaussian dynamics with a fixed basis}\label{sec:dynamics}
\begin{figure}[tbp]
\centering\includegraphics[width=\linewidth]{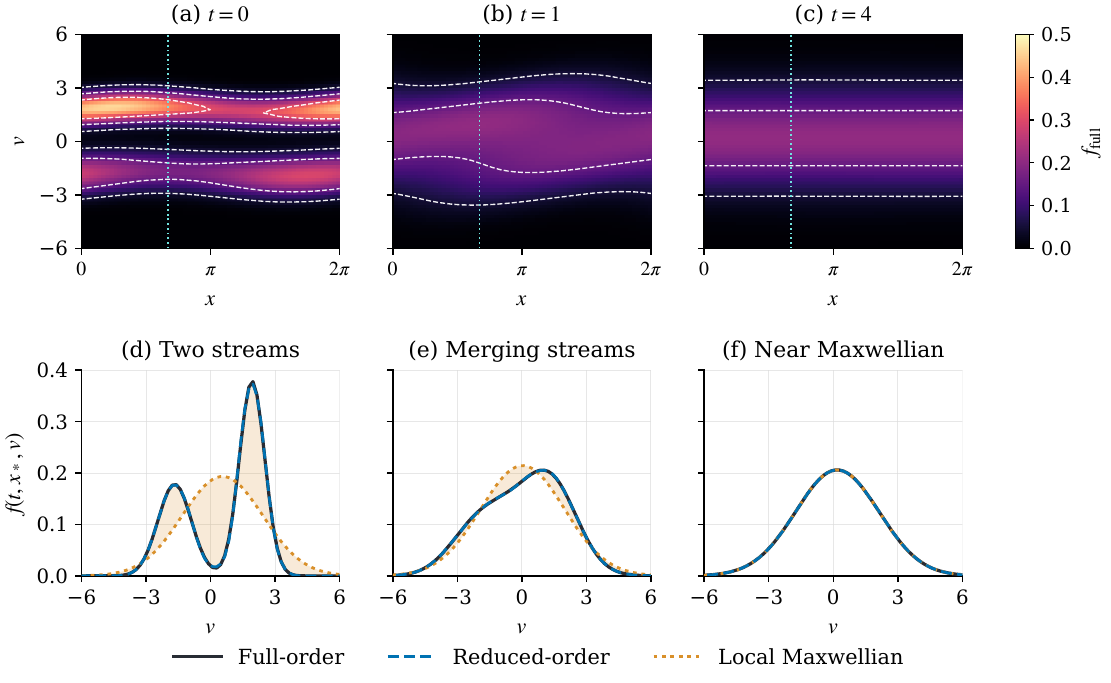}
\caption{Two-stream relaxation with the fixed $r=12$ post-trained basis
and positive ECM, $N_x=64$, $N_v=128$, $\Delta t=0.00125$.
(a)--(c) Full-order phase-space density on a common color scale; white
dashed contours show the reduced-order distribution at
$f=0.05,0.15,0.30$.  Cyan dotted lines locate $x_*=2.111$.
(d)--(f) Velocity profiles at that cell, with the local Maxwellian
$f_G$ of the full-order state from~\eqref{eq:nongaussian}.
Shading lies between the full-order profile and $f_G$, highlighting
the non-Gaussian information retained by the reduced model.}
\label{fig:evolution}
\end{figure}

In relaxation, the initially bimodal distribution approaches a local
Maxwellian as collisions merge the two streams.  The reduced-order
model follows both regimes with the same fixed basis
(Figure~\ref{fig:evolution}).  Its distribution error is largest in
the initial reconstruction and becomes smaller as the distribution
relaxes.

\begin{figure}[!tbp]
\centering\includegraphics[width=\linewidth]{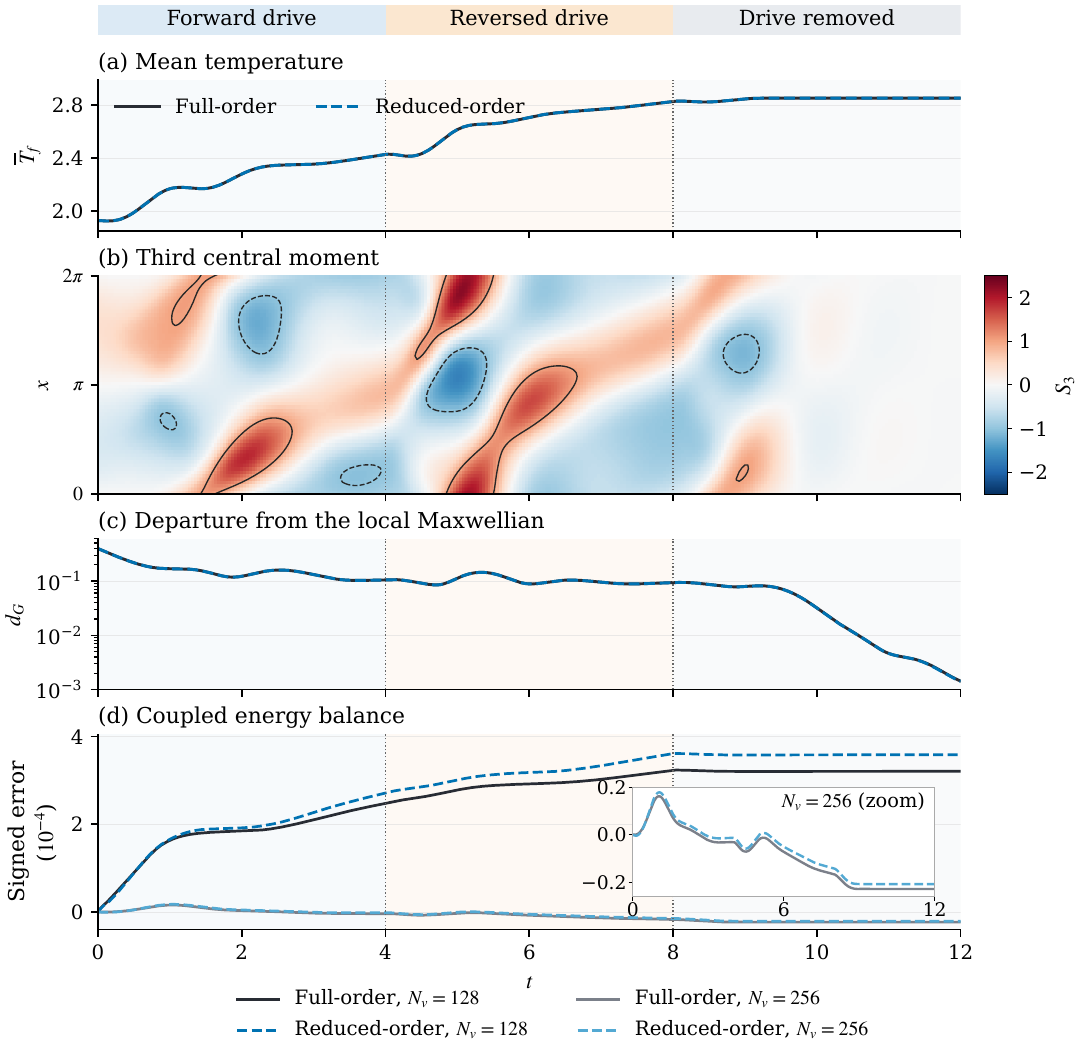}
\caption{Drive switching with the fixed $r=12$ post-trained basis
and positive ECM, $N_x=64$, $\Delta t=0.00125$.
The header marks the forcing stages of~\eqref{eq:switch-drive}.
(a) Mean temperature.  (b) Full-order third central moment $S_3(t,x)$;
black contours show the reduced-order values $S_3=\pm1$, with negative
values dashed.  (c) Deviation $d_G$ from the local Maxwellian.
(d) Signed energy balance error
$[\mathcal E(t)-\mathcal E(0)-W_{\rm ext}(t)]/\mathcal E(0)$,
cf.~\eqref{eq:defects}, on $N_v=128$ and $256$; the inset resolves the
fine-grid curves in the same vertical units.
Panels (a)--(c) use $N_v=128$.}
\label{fig:switch}
\end{figure}

Sustained heating raises the mean temperature and maintains a
non-Gaussian distribution beyond the training interval.  Drive
switching adds a change in the forcing history: reversal renews the
departure from the local Maxwellian, and removal allows relaxation
again (Figure~\ref{fig:switch}a--c).  The fixed reduced model tracks
these forced and relaxing phases without changing its basis.

To assess the number of data-adapted velocity functions, we compare
$r=8$ and $r=12$ network fits before post-training, using exact
collision sums.  The $r=8$ network is fitted separately to the first
six POD functions.  The Gaussian model ($r=2$), which retains only the
physical moments, provides a low-order reference.
For a same-dimensional analytic control in relaxation and heating,
we also use $r=12$ with
\[
 \psi_{2k-1}(v)=\sin(k\pi v/4),\qquad
 \psi_{2k}(v)=\cos(k\pi v/4),\qquad k=1,\ldots,5.
\]
All reconstructions start with the same density, mean velocity and
temperature.

\begin{table}[tbp]
\centering\small
\newcommand{\AnalyticRelaxationRow}{ & 12 (trig.) & 1.98 & 1.02 & 8.93 & 1.47 & 1.21 \\}
\newcommand{\AnalyticHeatingRow}{ & 12 (trig.) & 0.980 & 0.361 & 1.86 & 0.439 & 0.274 \\}
\caption{Relative errors, in percent, of the reduced-order model with
$r=2,8,12$ against the full-order model, $N_x=64$, $N_v=128$,
$\Delta t=0.00125$, exact sums.  The $r=2$ rows give the low-order
Gaussian reference.  `Initial' compares $f$ at $t=0$, the columns $f$, $J$ and $T_f$ give
the errors~\eqref{eq:error} over the whole run, and `late' uses $[2,4]$ for relaxation, $[6,8]$ for heating and $[8,12]$ for drive
switching.  Rows marked `trig.' use the fixed trigonometric basis defined in the text.}
\label{tab:rank}
\begin{tabular}{lrrrrrr}
\toprule
Problem & $r$ & initial $f$ & $f$ & $J$ & $T_f$ & late $f$\\
\midrule
Relaxation & 2 & 66.8 & 23.8 & 87.2 & 13.6 & 14.7 \\
 & 8 & 2.47 & 0.702 & 1.13 & 0.236 & 0.236 \\
 & 12 & 1.44 & 0.308 & 0.163 & 0.0398 & 0.0602 \\
\AnalyticRelaxationRow
\midrule
Heating & 2 & 39.6 & 22.8 & 85.7 & 19.1 & 23.6 \\
 & 8 & 1.50 & 1.38 & 2.34 & 0.377 & 0.617 \\
 & 12 & 0.799 & 0.340 & 0.271 & 0.0338 & 0.142 \\
\AnalyticHeatingRow
\midrule
Drive switching & 2 & 39.6 & 22.3 & 97.8 & 21.6 & 20.5 \\
 & 8 & 1.50 & 1.37 & 2.43 & 0.423 & 0.668 \\
 & 12 & 0.799 & 0.335 & 0.260 & 0.0411 & 0.144 \\
\bottomrule
\end{tabular}
\end{table}

Increasing $r$ from $8$ to $12$ improves the distribution and moment
predictions in all three problems (Table~\ref{tab:rank}).
At fixed $r=12$, the data-adapted functions also give smaller current
and temperature errors than the tested trigonometric control.
Section~\ref{sec:basis-study} compares the basis constructions at the
same dimension.

\subsection{Collision structure and global balance equations}\label{sec:structure}
The local collision invariants and the collision entropy identity are
preserved with each of the three bases and with positive quadrature.
Table~\ref{tab:experimental-structure} separates these local properties
from the balance errors of the coupled evolution.  All tested models
retain positive distributions, and collision entropy is nonincreasing
to solver precision.

Since $\int E\rho\,dx=0$ and $\frac{d}{dt}\frac12\int E^2\,dx=-\int EJ\,dx$
by~\eqref{eq:poisson} and the continuity equation,
integrating~\eqref{eq:moment-system} over $x$ shows that the total
momentum and energy change only through the external field:
\begin{equation}\label{eq:energy}
 \frac{d}{dt}\iint vf\,dx\,dv=\int E_{\rm ext}\rho\,dx,\qquad
 \frac{d}{dt}\Bigl(\frac12\iint v^2f\,dx\,dv+\frac12\int E^2\,dx\Bigr)
 =\int E_{\rm ext}J\,dx .
\end{equation}
The total mass, momentum and energy of the scheme are
\[
 \mathcal M=\Delta x\Delta v\sum_{i,j}f_{i,j},\qquad
 \mathcal J=\Delta x\Delta v\sum_{i,j}v_jf_{i,j},\qquad
 \mathcal E=\frac{\Delta x\Delta v}2\sum_{i,j}v_j^2f_{i,j}
 +\frac{\Delta x}2\sum_iE_i^2 .
\]
The external impulse $I_{\rm ext}$ and work $W_{\rm ext}$ accumulate
$\Delta x\sum_iE_{\rm ext}(t,x_i)\rho_i$ and
$\Delta x\sum_iE_{\rm ext}(t,x_i)J_i$ over the transport stages
of~\eqref{eq:splitting}, with their weights and with $A_{\rm ext}(t_n)$ in step $n$, and the balance errors
of the discrete counterparts
of~\eqref{eq:energy} are
\begin{equation}\label{eq:defects}
 \delta_{\mathcal J}(t)=\bigl|\mathcal J(t)-\mathcal J(0)-I_{\rm ext}(t)\bigr|,
 \qquad
 \delta_{\mathcal E}(t)=\frac{|\mathcal E(t)-\mathcal E(0)-W_{\rm ext}(t)|}{\mathcal E(0)} .
\end{equation}
Collisions keep every $\rho_i$, hence $E_i$, and conserve the local
invariants up to the stopping tolerance, so the balance errors change only in
transport, apart from this tolerance and round-off.  Transport conserves
$\mathcal M$, and since $\sum_iE_i\rho_i=0$, its momentum balance error comes from the
limited slopes in $v$ and the boundary values at $v=\pm12$. Its energy balance error also contains the velocity discretization of the work and of the
field energy.

\begin{table}[!htbp]
\centering\small
\setlength{\tabcolsep}{4pt}
\caption{Structure diagnostics on $N_x=64$, $N_v=256$ and $\Delta t=0.00125$.
All rows take maxima over the three two-stream tests.
The $r=12$ models use positive ECM.
The local collision conservation error $\epsilon_C^{\rm col}$ is the maximum
of $\|C^\top(q^{n+1}-q^n)\|_\infty$ over cells and collision steps, and $\epsilon_H^{\rm col}$ is the absolute
residual of the corresponding collision entropy identity.
The coupled diagnostics are $\max_t|\mathcal M(t)-\mathcal M(0)|/\mathcal M(0)$
and the maxima of~\eqref{eq:defects}.}
\label{tab:experimental-structure}
\begin{tabular}{lrrrrr}
\toprule
& \multicolumn{2}{c}{Collision step} & \multicolumn{3}{c}{Coupled evolution} \\
\cmidrule(lr){2-3}\cmidrule(lr){4-6}
Model & $\epsilon_C^{\rm col}$ & $\epsilon_H^{\rm col}$ & relative mass error & $\delta_{\mathcal J}$ & $\delta_{\mathcal E}$ \\
\midrule
Full-order & $2.47\times10^{-14}$ & $2.74\times10^{-14}$ & $2.99\times10^{-12}$ & $4.01\times10^{-5}$ & $2.29\times10^{-5}$ \\
POD & $5.87\times10^{-14}$ & $1.23\times10^{-13}$ & $7.49\times10^{-15}$ & $3.62\times10^{-5}$ & $2.28\times10^{-5}$ \\
Network fit & $5.96\times10^{-14}$ & $1.93\times10^{-13}$ & $5.94\times10^{-15}$ & $4.90\times10^{-5}$ & $2.08\times10^{-5}$ \\
Post-trained & $5.88\times10^{-14}$ & $1.16\times10^{-13}$ & $1.00\times10^{-14}$ & $4.00\times10^{-5}$ & $2.22\times10^{-5}$ \\
\bottomrule
\end{tabular}
\end{table}

The collision identities hold to solver tolerance, while the coupled
momentum and energy balance errors are comparable to those of the
full-order model, which uses the same transport discretization.
At the saved times, the total entropy
$\Delta x\Delta v\sum_{i,j}f_{i,j}\log f_{i,j}$ is nonincreasing
for all three bases.  Collision entropy dissipation is guaranteed by
Proposition~\ref{prop:sparse-collision}; this observed behavior of the
coupled method also includes the entropy error of MUSCL transport.

\subsection{Reduction error and discretization error}\label{sec:refinement}
To assess splitting and time discretization, we compare both split models with
an unsplit reference on $N_x=32$, $N_v=64$ up to $t=4$: the adaptive
eighth-order Dormand--Prince method DOP853
($\mathrm{rtol}=2\times10^{-11}$, $\mathrm{atol}=2\times10^{-12}$, maximum
step $10^{-3}$), applied in the variable $\log f$ to the semi-discrete
system $\frac{df_i}{dt}=R(f,t)_i-(\rho_i/\Delta v)\,\nu T_iK_{q_i}\log q_i$ that
the split full-order model discretizes in time.
Table~\ref{tab:experimental-refinement}a shows that the splitting and
time errors are much smaller than the reduction errors in both
relaxation and heating.

To measure the discretization errors, three hierarchies refine one
parameter at a time up to $t=4$: $N_x\in\{32,64,128\}$,
$N_v\in\{64,128,256\}$ and
$\Delta t\in\{0.0025,0.00125,0.000625,0.0003125\}$, with the other
parameters at their default values and the same network on every
velocity grid.  The observed order is the base-$2$ logarithm of the ratio
of the last two absolute space--time differences between successive
levels, with fine-grid values averaged to the coarser grid.
Table~\ref{tab:experimental-refinement}b shows approximately second-order
convergence in $x$ and $v$ and first-order convergence in time, for both
the distribution and the current.

\begin{table}[!htbp]
\centering\small
\setlength{\tabcolsep}{4pt}
\caption{Separating reduction from discretization error, up to $t=4$.
(a) Errors~\eqref{eq:error} against the unsplit DOP853 reference on
$N_x=32$, $N_v=64$, with $\Delta t=0.00125$ for the split models.
(b) Observed orders $p_f$ and $p_J$ and last absolute space--time
differences $\Delta_f$, with F/R denoting full-order/reduced.
The reduced-model error $e_f$ uses the full-order model at the
finest level of each hierarchy.  Both models use the same grids and
time steps.  The reduced model uses the network approximation and
exact collision sums throughout.}
\label{tab:experimental-refinement}
\begin{tabular}{lrrrr}
\toprule
\multicolumn{5}{l}{(a) Space--time errors against the unsplit reference (\%)} \\
& \multicolumn{2}{c}{$f$} & \multicolumn{2}{c}{$J$} \\
\cmidrule(lr){2-3}\cmidrule(lr){4-5}
Problem & Full-order & Reduced & Full-order & Reduced \\
\midrule
Relaxation & 0.0167 & 0.3000 & 0.0123 & 0.1518 \\
Heating & 0.0056 & 0.4487 & 0.0106 & 0.2651 \\
\bottomrule
\end{tabular}
\par\medskip
\begin{tabular}{llrrrrr}
\toprule
\multicolumn{7}{l}{(b) Orders, last absolute differences and reduction error} \\
Problem & Refine & $p_f$ (F/R) & $p_J$ (F/R) & $\Delta_f^{\rm F}$ & $\Delta_f^{\rm R}$ & $e_f$ (\%) \\
\midrule
Relaxation & $x$ & 2.02 / 2.02 & 1.85 / 1.93 & $2.06\times10^{-4}$ & $2.09\times10^{-4}$ & 0.308 \\
 & $v$ & 1.99 / 1.98 & 1.89 / 1.95 & $2.08\times10^{-3}$ & $2.09\times10^{-3}$ & 0.311 \\
 & $t$ & 1.00 / 1.00 & 0.96 / 0.96 & $8.18\times10^{-5}$ & $8.28\times10^{-5}$ & 0.308 \\
\midrule
Heating & $x$ & 2.26 / 2.19 & 2.13 / 2.02 & $7.81\times10^{-4}$ & $8.23\times10^{-4}$ & 0.448 \\
 & $v$ & 1.86 / 1.85 & 1.97 / 2.04 & $2.69\times10^{-3}$ & $2.67\times10^{-3}$ & 0.435 \\
 & $t$ & 1.00 / 1.00 & 0.98 / 0.98 & $3.28\times10^{-5}$ & $3.46\times10^{-5}$ & 0.443 \\
\bottomrule
\end{tabular}
\end{table}
Against the finest full-order solution of each hierarchy, the
reduced-model error changes little under space or time refinement.
Velocity refinement removes the additional coarse-grid error and
approaches the same error level.  These refinements indicate that
reduction, including representation and evolution within the fixed
family, controls the remaining error.
The coupled energy balance is more sensitive to the velocity grid:
Figure~\ref{fig:switch}d shows its improvement on the finer grid through
all three forcing stages.  Further velocity refinement in the matched
cost study below changes $f$ by about $\MatchedFineGridDifference\%$
from $N_v=512$ to $1024$, in both the full-order model and the
post-trained model with positive quadrature.

The collision discretization is also checked against an analytic
homogeneous Dougherty solution in Appendix~\ref{app:implementation}.
The discrete collision right-hand side converges at second order in
velocity to the analytic generator, independently of time integration;
Table~\ref{tab:ou-analytic} reports distribution and moment errors for
the full-order and final reduced models.

\begin{remark}
Enlarging the velocity interval from $[-12,12]$ to $[-18,18]$ at fixed
$\Delta v$ changes the full-order solution on the common interval by
at most $2\times10^{-9}$ in relaxation and heating up to $t=4$.
At $\Delta v=24/256$, the same comparison over the complete heating
($t=8$) and drive-switching ($t=12$) trajectories gives a relative
space--time difference~\eqref{eq:error} below $\DomainFullLongDelta$.
\end{remark}

\subsection{Choice and training of the basis}\label{sec:basis-study}
We compare the three bases on $N_v=128$, $256$ and $512$, all with
$r=12$ and positive quadrature.  The ECM rules use a selection budget of
$64$ per rule, retaining $63$ or $64$ nodes and $62$ to $64$ interfaces on
every grid. Across the nine basis--grid pairs,
the largest relative error of the collision matrices on the validation
probability vectors is $\DynamicRuleError$ in the Frobenius norm.
A fixed number of quadrature points keeps the accuracy on the finer grids:
for each basis and problem, every error in
Table~\ref{tab:dynamic-accuracy} changes by at most $3.1\%$ of its value
between $N_v=128$ and $512$.
For the final post-trained basis, ECM alters the solution by much less
than the reduction error (Table~\ref{tab:final-ecm}).

Table~\ref{tab:dynamic-accuracy} compares the three bases on $N_v=256$.
Post-training lowers the space--time distribution error of the network
approximation in all three problems.  It also improves on POD in
heating and drive switching, while POD gives the smallest distribution
error in relaxation.  The moment errors show a problem-dependent
ordering.

\begin{figure}[tb]
\centering\includegraphics[width=\linewidth]{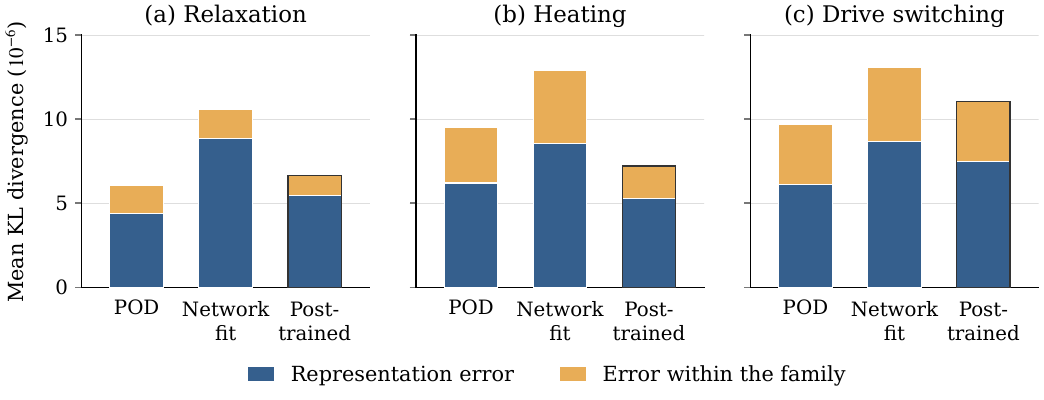}
\caption{Error decomposition for three fixed bases, $r=12$, $N_x=64$,
$N_v=256$ and $\Delta t=0.00125$, with positive collision quadrature
and full-grid MUSCL.  Each bar is the mean KL error of the normalized
cell probability vectors.  Its two parts are the representation error
and the error within the family in~\eqref{eq:kl-error-decomposition},
averaged over all cells and saved times with trapezoidal time weights.}
\label{fig:basis-training}
\end{figure}

Figure~\ref{fig:basis-training} explains the change through the KL
decomposition~\eqref{eq:kl-error-decomposition}.  We average each term
over space and time,
\[
 \overline{\KL}(a\Vert b)
 =\frac{1}{T_{\rm end}N_x}\int_0^{T_{\rm end}}
   \sum_i\KL\bigl(a_i(t)\Vert b_i(t)\bigr)\,dt,
\]
where $a_i$ and $b_i$ are cell probability vectors and $T_{\rm end}$ is
the final time of the test problem.  Relative to the network fit,
post-training reduces both the representation error and the error
within the family in all three problems.  The representation error
remains the larger component.  Thus training through the dynamics
improves both which distributions the basis can represent and how
the reduced trajectory follows their best reconstruction in these
tests.

\begin{table}[tb]
\centering\small
\setlength{\tabcolsep}{4pt}
\caption{Three fixed bases with positive quadrature, $r=12$,
$N_x=64$, $N_v=256$, $\Delta t=0.00125$: space--time relative
errors~\eqref{eq:error}, in percent, against the full-order model.  The $f$,
$J$, $T_f$ and $S_3$ columns use the whole run, the last column its last
quarter: $[3,4]$, $[6,8]$ and $[9,12]$.}
\label{tab:dynamic-accuracy}
\begin{tabular}{llrrrrr}
\toprule
Problem & Basis & $f$ & $J$ & $T_f$ & $S_3$ & last-quarter $f$ \\
\midrule
Relaxation & POD & 0.193 & 0.212 & 0.049 & 0.642 & 0.035 \\
 & Network fit & 0.311 & 0.166 & 0.040 & 0.693 & 0.028 \\
 & Post-trained & 0.217 & 0.176 & 0.038 & 0.469 & 0.040 \\
\midrule
Heating & POD & 0.283 & 0.169 & 0.030 & 0.329 & 0.145 \\
 & Network fit & 0.335 & 0.272 & 0.034 & 0.548 & 0.140 \\
 & Post-trained & 0.210 & 0.167 & 0.037 & 0.368 & 0.097 \\
\midrule
Drive switching & POD & 0.276 & 0.194 & 0.039 & 0.501 & 0.115 \\
 & Network fit & 0.329 & 0.261 & 0.041 & 0.589 & 0.122 \\
 & Post-trained & 0.241 & 0.264 & 0.063 & 0.634 & 0.136 \\
\bottomrule
\end{tabular}
\end{table}

\subsection{Online cost}\label{sec:cost}
We time the full-order model and the reduced-order model with positive quadrature
and the network approximation or the post-trained basis on drive switching up to
$t=12$, with $N_x=64$, $\Delta t=0.00125$ and $N_v=128,256,512,1024$. The
reduced computations have $r=12$, at most $64$ positive nodes and interfaces,
and the same full-grid MUSCL transport as the full-order model.  The full-order model solves each collision step by a damped Newton method whose matrix
$\diag(e^y)+h\nu T_{q^n}K_{q^n}$ in $y=\log q^{n+1}$ is, by~\eqref{eq:mobility},
a tridiagonal matrix minus a matrix of rank at most two, so that each
iteration costs $O(N_v)$ operations.  Each configuration is repeated three times in one computing job on an AMD
EPYC 9654 node, with one BLAS thread per process. Each repetition stays
on one physical core and rotates the grid and method order.  Online time
includes the initial reconstruction, time stepping and the diagnostics
common to both models. Setup, three warm-up steps and output I/O are
excluded.

\begin{figure}[!htbp]
\centering\includegraphics[width=\linewidth]{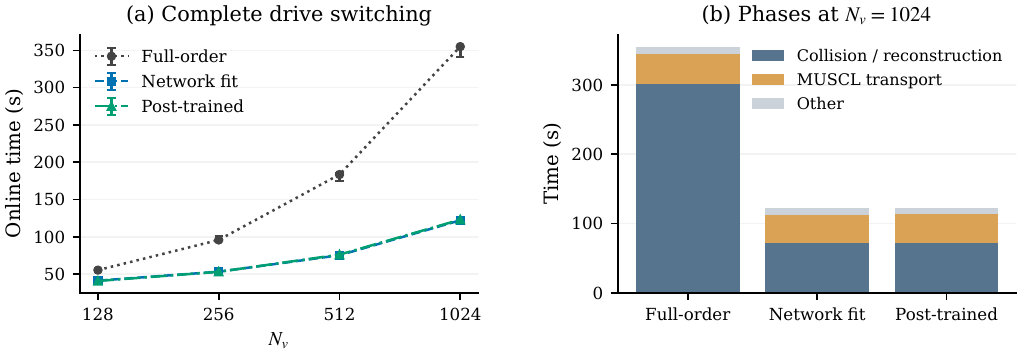}
\caption{Online costs for drive switching.  (a) Median time and
minimum--maximum range of three runs.  (b) Collision, MUSCL on the full velocity grid
and remaining costs of the median-time run for each model at
$N_v=1024$.  Collision time includes reconstruction and its diagnostics.
The remainder includes the global diagnostics.
The post-trained model is $\MatchedDynamicSpeedMinPrecise$--$\MatchedDynamicSpeedMaxPrecise$
times faster in median time.  At $N_v=1024$, the distribution errors
against the same-grid full-order model are $\MatchedDynamicFineF\%$
(post-trained) and $\MatchedNeuralFineF\%$ (network fit).}
\label{fig:matched-cost}
\end{figure}

\begin{samepage}
Post-training improves distribution accuracy at comparable online cost
(Figure~\ref{fig:matched-cost}).  ECM lowers collision time relative
to optimized exact sums for the same final basis
(Table~\ref{tab:final-ecm}), with the largest online gain on the finest
grid.

\end{samepage}

\FloatBarrier
\section{Conclusions}\label{sec:conclusions}
We have developed a low-dimensional entropy-based moment method for
non-Gaussian Vlasov--Poisson--Dougherty dynamics.  The model augments
density, mean velocity and temperature with data-adapted velocity
moments.  For every fixed basis satisfying the stated assumptions,
the reduced collision step is uniquely solvable, positive,
conservative in mass, momentum and kinetic energy, and
entropy-dissipative.  The variational formulation gives a
best-approximation characterization relative to the full collision
step with the same input and frozen mobility.

The same fixed basis tracks the distribution and macroscopic responses
in two-stream relaxation, sustained heating and drive switching.  In these tests,
dynamics post-training lowers the distribution error of the network
fit and gives smaller distribution errors than POD in the two driven
tests.  The mean KL decomposition also shows
improvements in both representation and evolution within the family
relative to the network fit.

Positive empirical cubature preserves the collision invariants and
entropy dissipation while assembling the collision matrix and the
approximate Newton matrix from small quadrature rules.  With
full-grid MUSCL retained, matched drive-switching timings show lower
online cost than the full-order solver.  The method thus separates
basis design from the collision guarantees.  Extending it to several
velocity dimensions and reducing the remaining full-grid operations
are natural next steps.

\section*{Acknowledgment}
Guang Lin gratefully acknowledges the support from DARPA MATHBAC program under grant number: HR0011-26-3-E152, the National Science Foundation (DMS-2533878, DMS-2053746, DMS-2134209, ECCS-2328241, CBET-2347401 and OAC-2311848), and U.S.~Department of Energy (DOE) Office of Science Advanced Scientific Computing Research program under the "Uncertainty Quantification for Multifidelity Operator Learning (MOLUcQ)" project (Project No. 81739), DE-SC0023161, the SciDAC LEADS Institute, and DOE–Fusion Energy Science, under grant number: DE-SC0024583.

\section*{Code and data availability}
The implementation, trained bases, cubature rules and supporting
numerical data are available from the authors upon request.

\appendix
\makeatletter
\renewcommand\@seccntformat[1]{%
 \@ifundefined{#1@cntformat}{\csname the#1\endcsname\quad}%
  {\csname #1@cntformat\endcsname}}
\newcommand\section@cntformat{\appendixname~\thesection\quad}
\makeatother

\section{Training details}\label{app:training}

\paragraph{Data.}
The training and validation trajectories use $N_x=32$, $N_v=128$ and
$\Delta t=0.0025$ up to $t=4$.
The $6{,}464$ POD probability vectors come from the training trajectories at every
saved time (spacing $0.04$) and every second cell.  The $176$ validation
probability vectors come from the validation trajectories at every tenth saved time
and every fourth cell.  The rules of Section~\ref{sec:rules} are fitted
to $382$ probability vectors: $352$ of the training trajectories at every tenth saved
time and every fourth cell, interpolated in $\log q$ to other grids, and
$30$ Gaussians with means $-2,-1,0,1,2$ and temperatures $0.25,0.5,1,2,4,6$,
evaluated on each grid.  After exponentiation, every interpolated
profile is renormalized on the target grid.  Nodes outside the source
cell-center range use cubic extrapolation of the endpoint spline segments.

\paragraph{Network approximation of the POD basis.}
The network $\mathcal N_\vartheta$ takes input $v/4$, has two hidden
$\tanh$ layers of width $64$, and has $d$ outputs.
The fixed output scales and fitting weights in~\eqref{eq:neural-pod} are
\[
 \gamma_k^2=\frac1{N_v}\sum_j(\Psi_{\rm POD})_{jk}^2,\qquad
 w_j=\frac1{2N_v}+\frac{\omega_j}2.
\]
The parameters are initialized by the default rule of PyTorch with
random seed $619$.  Training is in double precision and runs $18{,}000$
full-batch Adam updates on~\eqref{eq:neural-pod} with learning rate $0.002$, multiplied by $0.3$
after updates $5{,}001$ and $9{,}001$.  Networks are saved every
$3{,}000$ updates.  The validation rule of Section~\ref{sec:neural-pod}
requires successful moment matching on all $176$ vectors and selects
update $18{,}000$ for $d=10$ and $12{,}000$ for $d=6$
(Section~\ref{sec:dynamics}).

\paragraph{Post-training.}
The initial reconstruction~\eqref{eq:moment-matching} is differentiated
with the Jacobian~\eqref{eq:moment-jacobian}.  Gradients pass through
the coupled steps, weighted QR, the Poisson solve and the active
branches of the MUSCL limiter; the basis rank is checked throughout.
For a converged root $\widetilde F_\vartheta(z;q^n)=0$ with
$h=\Delta t/2$, the adjoint variable $\mu\in\R^r$ satisfies
\begin{equation}\label{eq:implicit-training-gradient}
 \begin{aligned}
 \partial_z\widetilde F_\vartheta
 &=\nabla m_\vartheta(z)+h\nu T_{q^n}\widetilde L_{q^n,\vartheta},\\
 (\partial_z\widetilde F_\vartheta)^\top\mu
 &=\nabla_z\mathcal L_{\rm dyn},\\
 (\nabla_z\mathcal L_{\rm dyn})^\top\,\mathrm dz
 &=-\mu^\top\bigl(
 \partial_\vartheta\widetilde F_\vartheta\,\mathrm d\vartheta
 +\partial_{q^n}\widetilde F_\vartheta\,\mathrm dq^n\bigr).
 \end{aligned}
\end{equation}
Here $\nabla_z\mathcal L_{\rm dyn}$ is the loss gradient propagated to
this root from subsequent computations, and partial derivatives hold
the other root arguments fixed.  The chain rule combines these terms
over all roots with the explicit dependence on $B_\vartheta$.

Training starts from the $d=10$ network fitted to the POD basis.
For its output weights and biases, we optimize
\[
 W_{k,:}=s_k\widehat W_{k,:},\qquad b_k=s_k\widehat b_k,\qquad
 s_k^2=\frac1{64}\|W_{k,:}^{(0)}\|_2^2+(b_k^{(0)})^2,
\]
with $s_k$ fixed at the initial checkpoint.
The weights and biases of all three affine layers are trainable;
the physical functions $1,v,v^2$ and the scales $s_k$ are fixed.
The $200$
training intervals start at $t=0,0.08,\ldots,3.92$ on each training
trajectory and the twelve validation intervals at
$t=0,0.4,1.2,2.4,3.2,3.92$ on each validation trajectory.  The mean squares
in~\eqref{eq:dynamic-loss} are $\overline{J^2}=0.0800407$ and
$\overline{S_3^2}=1.3999433$.  We use double precision and $600$ Adam iterations
with the PyTorch defaults, one interval per update, with a shuffled
traversal of all $200$ intervals on each pass and random seed $20261028$.
The initial learning rate is $10^{-4}$, multiplied by $0.3$ after updates
$300$ and $450$, and the gradient norm is capped at $1$.  An update is
halved, at most five times, until its interval objective is at most $1.25$
times the previous one and all nonlinear solves converge; otherwise,
the parameters and optimizer state are restored.  Of the
networks after $0$, $50$, $100$, $200$, $300$, $450$ and $600$ updates,
the best three on the validation intervals with the network approximation rules are
evaluated again with rules refitted to the $382$ probability vectors, and the best is kept.  This selects iteration $600$.
The implicit
derivatives~\eqref{eq:implicit-training-gradient} agree with central
differences to relative errors of $1.3\times10^{-9}$ for one collision
step and $1.0\times10^{-6}$ through three coupled steps.

\paragraph{KL diagnostics.}
For every saved reference state in Figure~\ref{fig:basis-training}, we
solve~\eqref{eq:moment-matching} with the corresponding fixed basis on
the full velocity grid, using the saved reduced coefficients as the
initial guess and the reconstruction tolerances of
Appendix~\ref{app:implementation}.  All cells and saved times are used;
the two terms of~\eqref{eq:kl-error-decomposition} are evaluated from
these reconstructions and the saved reduced states.

\section{Implementation details}\label{app:implementation}

\paragraph{Basis coordinates.}
On the training grid $N_v=128$, let $\Psi_{jk}=\psi_k(v_j)$ and use
\begin{equation}\label{eq:basis-qr}
 \Omega^{1/2}[\one,v,v^2,\Psi]=\widehat Q\widehat R,\qquad
 B=[\one,v,v^2,\Psi]\,(\widehat R^{-1})_{:,2:r+1}.
\end{equation}
The weights in $\Omega=\diag(\omega_j)$ are
$\omega_j=e^{-v_j^2/8}/\sum_s e^{-v_s^2/8}$.
Then $B^\top\Omega B=I_r$.  The factor $\widehat R$ is reused when
evaluating the functions on other grids.
The ordering of the QR columns gives
\begin{equation}\label{eq:qr-physical-columns}
 B_{:,1}\in\spn\{\one,v\},\qquad
 \spn\{\one,B_{:,1},B_{:,2}\}=\range(C).
\end{equation}
\paragraph{Newton solvers.}
Newton's method for~\eqref{eq:full} acts on $y=\log q^{n+1}$ with the
residual
\[
 F_{\rm full}(y)=e^y-q^n+h\nu T_{q^n}K_{q^n}y.
\]
For the Newton correction $\delta y$, backtracking takes the first
$\alpha=2^{-k}$, $k=0,\ldots,34$, such that
\[
 \|F_{\rm full}(y+\alpha\delta y)\|_2
 <(1-10^{-4}\alpha)\|F_{\rm full}(y)\|_2.
\]
The same test is used for $F$ and $\widetilde F$ with corrections
$\delta z$.  Trials with $\min_j y_j<-700$ in the full solve or
$\min_j\log p_j\le-700$ in a reduced solve are rejected.
If all $35$ trials fail in~\eqref{eq:sparse-newton}, the exact
covariance~\eqref{eq:moment-jacobian} replaces the sampled covariance
for the unaccepted cells, and the same $35$ trials are repeated.
The frozen collision matrix is unchanged.  Failure of this retry
terminates the run with an error; the production exact-sum solves
terminate if their search fails.  The outer loops have $45$ iterations for the
full solve and $100$ for the reduced solve, each beginning with the
stopping test; exhausting either loop also terminates the run.
The full solve stops when $\|F_{\rm full}\|_\infty\le4\times10^{-14}$ and
\[
 \max\Bigl\{|\one^\top(q^{n+1}-q^n)|,\ \frac{|v^\top(q^{n+1}-q^n)|}{1+|v^\top q^n|},\
 \frac{|(v^2/2)^\top(q^{n+1}-q^n)|}{1+(v^2/2)^\top q^n}\Bigr\}\le8\times10^{-15}.
\]
The reduced solves, with exact sums, with quadrature and in training,
stop when the residual is at most $4\times10^{-13}$ in the maximum norm
and its first two components, which equal $B_{:,1:2}^\top(p(z)-q^n)$
by~\eqref{eq:qr-physical-columns} and $K_qC=\widetilde K_qC=0$, are at most
$2\times10^{-14}$.
The exact covariance~\eqref{eq:moment-jacobian} is used for the
reconstruction~\eqref{eq:moment-matching} and as the line-search fallback
for~\eqref{eq:sparse-newton}.

\paragraph{Empirical cubature rules.}
For each fixed basis and velocity grid, the node rule is fitted to the
normalization $\sum_jq_j$, the moments $B^\top q$ and the distinct entries
of $\sum_jq_jB_j^\top B_j$ over the training probability vectors.
The interface rule fits the distinct entries of
\[
 \sum_e\Lambda(q_e,q_{e+1})\,[G,DB]_{e,:}^\top[G,DB]_{e,:},
\]
which contain $G^\top W_qG$, $G^\top W_qDB$ and $(DB)^\top W_qDB$.
Each fitted entry is divided by its root-mean-square value over the
probability vectors, or by $10^{-2}$ times the largest such value if that
is larger.  The scaled training matrix is compressed in its row space
by an SVD, retaining singular values larger than $10^{-13}$ times the
largest.  Each of the $64$ greedy iterations selects the remaining
column with largest normalized correlation to the current residual and
recomputes the selected weights by NNLS.  Nodes or interfaces with
weight at most $10^{-12}$ are removed after the greedy steps.
The reduced collision equilibria are unchanged when the interfaces satisfy
\begin{equation}\label{eq:edge-rank}
 \rank[G_\Gamma,D_\Gamma B]-\rank G_\Gamma=r-2.
\end{equation}
Since $\ker(I-P_\Gamma)=\range(G_\Gamma)$, this gives
\[
 \ker\widetilde L_q=\ker L_q=\{z:Bz\in\range(C)\}.
\]

The retained rules are checked for linear independence on $\mathcal N$,
$\rank G_\Gamma=2$ and~\eqref{eq:edge-rank}.
For the post-training basis and rule checks, the numerical rank of an
$m\times n$ matrix counts singular values satisfying
$\sigma_k>\max(m,n)\epsilon_{\rm mach}\sigma_1$, with
$\epsilon_{\rm mach}$ the double-precision machine epsilon.

\paragraph{Final-basis cubature check.}
In Table~\ref{tab:final-ecm}, $\delta_g$ uses~\eqref{eq:error} with
$g_{\rm red}=g_{\rm ECM}$ and $g_{\rm full}=g_{\rm exact}$,
for $g\in\{f,J,S_3\}$.  Both implementations use cached matrix
products and the same stopping tests.  Across the final-basis runs,
collision invariant and entropy-identity defects are below
$\FinalECMInvariant$ and $\FinalECMEntropyIdentity$, respectively.
No exact-covariance fallback was triggered.
In the timing runs, the mean of the maximum Newton count over cells
per batched call ranges from $\FinalECMIterationMin$ to
$\FinalECMIterationMax$; the initial reconstruction is included.
Global impulse and work balances use~\eqref{eq:defects}.

\begin{table}[!ht]
\centering\small
\setlength{\tabcolsep}{4pt}
\caption{ECM versus exact sums for the final post-trained basis, $r=12$.
Accuracy columns give dimensionless space--time relative differences,
maximized over
the three problems at each grid.  Times (exact / ECM, seconds) are
medians of three drive-switch runs on one AMD EPYC 9654 node,
$N_x=64$, $\Delta t=0.00125$.  Collision includes reconstruction.}
\label{tab:final-ecm}
\begin{tabular}{rccc cc}
\toprule
$N_v$ & $\max\delta_f$ & $\max\delta_J$ & $\max\delta_{S_3}$ & Online (s) & Collision (s)\\
 & & & & exact / ECM & exact / ECM\\
\midrule
128 & $7.4\times10^{-7}$ & $1.7\times10^{-6}$ & $3.9\times10^{-6}$ & 46.8 / 42.5 & 36.2 / 32.3\\
256 & $1.2\times10^{-6}$ & $1.3\times10^{-6}$ & $3.4\times10^{-6}$ & 69.2 / 58.3 & 52.6 / 41.2\\
512 & $1.8\times10^{-6}$ & $1.9\times10^{-6}$ & $4.4\times10^{-6}$ & 120.9 / 83.2 & 90.9 / 53.2\\
\bottomrule
\end{tabular}
\end{table}

\paragraph{Analytic collision benchmark.}
For homogeneous Dougherty relaxation, $\rho=1$, $u$ and $T$ are constant.
A Gaussian mixture has the exact whole-line solution
\[
 \begin{aligned}
 f(t,v)&=\sum_{\ell=1}^{2}w_\ell M_{1,m_\ell(t),\theta_\ell(t)}(v),\\
 m_\ell(t)&=u+(m_\ell(0)-u)e^{-\nu t},\\
 \theta_\ell(t)&=T+(\theta_\ell(0)-T)e^{-2\nu t}.
 \end{aligned}
\]
We take $(w_1,w_2)=(0.6,0.4)$, $(m_1(0),m_2(0))=(1.8,-2.1)$,
$(\theta_1(0),\theta_2(0))=(0.45,0.7)$ and $\nu=0.4$;
$u=0.24$ and $T=4.2004$ follow from the initial moments.
The third central moment and fourth cumulant satisfy
\[
 S_3(t)=S_3(0)e^{-3\nu t},\qquad
 \kappa_4(t)=\kappa_4(0)e^{-4\nu t},\qquad
 \kappa_4=\int(v-u)^4f\,dv-3T^2.
\]
Table~\ref{tab:ou-analytic} compares the frozen-mobility full-order step
and the final $r=12$ basis, with exact sums or positive ECM, on
$[-12,12]$.  The right-hand-side test uses the analytic $\partial_tf$
at $t=0.2$ and $N_v=64,128,256,512$.

\begin{table}[!ht]
\centering\small
\caption{Homogeneous collision errors against the analytic solution,
$0\le t\le1$, $\Delta t=0.0005$.  All entries are relative space--time
errors in percent, with time samples spaced by $0.02$.
Both ROM variants use the same final post-trained basis; ECM uses the
frozen positive rules with a budget of $64$.}
\label{tab:ou-analytic}
\begin{tabular}{rlrrr}
\toprule
$N_v$ & Model & $e(f)$ (\%) & $e(S_3)$ (\%) & $e(\kappa_4)$ (\%)\\
\midrule
128 & Full-order & 0.15 & 0.11 & 0.0912 \\
 & ROM, exact sums & 0.277 & 0.0504 & 0.179 \\
 & ROM, ECM & 0.277 & 0.0504 & 0.179 \\
\midrule
256 & Full-order & 0.0549 & 0.0833 & 0.0736 \\
 & ROM, exact sums & 0.24 & 0.0261 & 0.16 \\
 & ROM, ECM & 0.24 & 0.0261 & 0.16 \\
\midrule
512 & Full-order & 0.0326 & 0.0769 & 0.0695 \\
 & ROM, exact sums & 0.237 & 0.0238 & 0.156 \\
 & ROM, ECM & 0.237 & 0.0238 & 0.156 \\
\bottomrule
\end{tabular}
\end{table}

\paragraph{Transport.}
The Poisson solve divides each discrete Fourier coefficient of $\rho$
with nonzero wavenumber below $N_x/2$ by $\mathrm i$ times the wavenumber
and sets the others to zero.  Each term of the discrete Parseval sum for
$\sum_iE_i\rho_i$ is then purely imaginary and odd in the wavenumber, so
$\sum_iE_i\rho_i=0$.  With $v_j^\pm$ and $a_i^\pm$ the positive and negative
parts of $v_j$ and $a_i$, the upwind fluxes in $R(f,t)$ are
\[
 \mathcal F^x_{i+1/2,j}=v_j^+\bigl(f_{i,j}+\tfrac12\sigma^x_{i,j}\bigr)
 +v_j^-\bigl(f_{i+1,j}-\tfrac12\sigma^x_{i+1,j}\bigr),\qquad
 \mathcal F^v_{i,j+1/2}=a_i^+\bigl(f_{i,j}+\tfrac12\sigma^v_{i,j}\bigr)
 +a_i^-\bigl(f_{i,j+1}-\tfrac12\sigma^v_{i,j+1}\bigr),
\]
with slopes
$\sigma=\operatorname{minmod}(1.5\,\Delta_-f,\tfrac12(\Delta_-f+\Delta_+f),1.5\,\Delta_+f)$
from the backward and forward differences in $i$ (periodic) for $\sigma^x$
and in $j$ for $\sigma^v$, and $\sigma^v=0$ at $j=1,N_v$.  The function
$\operatorname{minmod}$ returns the argument of smallest modulus if all
three have one sign and $0$ otherwise.  Both slopes entering the
difference of two upwind values lie between $0$ and $1.5$ times the
difference of the cell values, so the difference of the reconstructed
values is between $\frac14$ and $\frac74$ times that of the cell values,
and at $j=1,N_v$ the single inflow value is a convex combination of two
neighboring values.  Hence, under the bound stated
after~\eqref{eq:splitting}, $f_{i,j}+\Delta t\,R(f,t)_{i,j}$ has
nonnegative weights on neighboring values and a weight on $f_{i,j}$ of at least
\[
 1-\frac74\Delta t\left(\frac{|v_j|}{\Delta x}
                         +\frac{|a_i|}{\Delta v}\right)
 \ge1-\frac74\cdot0.48=0.16.
\]



\end{document}